\documentclass[final]{siamonline1116}
\usepackage{multirow}

\usepackage{lipsum}
\usepackage{array}
\usepackage{amsfonts}
\usepackage{graphicx}
\usepackage{epstopdf}
\usepackage{algorithmic}
\usepackage{subfigure}
\usepackage{bbm}

\graphicspath{{./}}

\ifpdf
  \DeclareGraphicsExtensions{.eps,.pdf,.png,.jpg}
\else
  \DeclareGraphicsExtensions{.eps}
\fi

\numberwithin{theorem}{section}

\newcommand{\TheTitle}{Statistical Treatment of Inverse Problems
Constrained by Differential Equations-Based Models with Stochastic
Terms}
\newcommand{\TheAuthors}{E.~M.~Constantinescu, 
      N.~Petra,
      J.~Bessac, 
      C.~G.~Petra}

\headers{Statistical Treatment of Inverse Problems with
Stochastic Terms}{\TheAuthors}

\title{{\TheTitle}\thanks{Submitted to the editors on October 13, 2018.
\funding{The work on the general scoring methodology, statistical analysis, computational experiments related to power grid, and design of the overall computational framework was supported by the U.S. Department of Energy, Office of Science, Advanced Scientific Computing Research Program under contracts DE-AC02-06CH11357 (at Argonne) and  DE-AC52-07NA27344 (at Lawrence Livermore). In addition, the National Science Foundation (NSF) grant SI2-SSI ACI-1550547 funded the hIPPYlib-related developments and the NSF grant  CAREER-1654311 supported the mathematical and computational developments as well as the computational experiments related to inversion governed by PDEs.
}}}

\author{
      Emil M. Constantinescu\footnotemark[2]\footnotetext[2]{Mathematics and Computer Science
Division, Argonne National Laboratory, Lemont, IL 60439
        (\email{emconsta@mcs.anl.gov}, \email{jbessac@anl.gov}).}
      \and
       No\'emi Petra\footnotemark[4]~\footnotemark[5]\footnotetext[4]{Applied Mathematics Department, University of
California, Merced, Merced, CA 95340
        (\email{npetra@ucmerced.edu}).}\footnotetext[5]{Corresponding author.}
        \and
      Julie Bessac\footnotemark[2]
      \and
      Cosmin G. Petra\footnotemark[3]\footnotetext[3]{Center for Applied Scientific Computing,
Lawrence Livermore National Laboratory, Livermore, CA 94550
(\email{petra1@llnl.gov}).}
}

\usepackage{amsopn}

\usepackage{booktabs}
\usepackage[mathscr]{euscript}
\usepackage{color}
\usepackage{graphicx}
\usepackage{algorithmic,algorithm}
\usepackage{tikz,pgfplots}
\usepackage{upgreek}
\usepackage{multirow}
\usepackage{yfonts}
\usepackage{mathtools}
\usepackage[normalem]{ulem}
\usepackage{latexsym}
\usepackage{url}
\usepackage{amsmath,amssymb,amsbsy,amsfonts}
\usepackage{enumitem}

\definecolor{mygreen}{rgb}{0.5,0.05,0.15}

\definecolor{utorange}{rgb}{0.8,0.33,0.}
\definecolor{themec}{RGB}{51,108,121}
\definecolor{darkred}{rgb}{.6,.1,.1}
\definecolor{darkblue}{rgb}{.1,.1,.9}
\definecolor{greenback}{rgb}{.19,.94,.13}
\definecolor{orange}{rgb}{.76,.39,.13}
\definecolor{grass}{rgb}{.19,.64,.13}
\definecolor{sierp}{RGB}{209,28,209}
\definecolor{bgorange}{rgb}{1.,.95,.78}
\definecolor{grassgreen}{RGB}{92,135,39}
\definecolor{thinbox}{rgb}{.7,.8,1.}

\renewcommand{\vec}[1]{{\mathchoice
                     {\mbox{\boldmath$\displaystyle{#1}$}}
                     {\mbox{\boldmath$\textstyle{#1}$}}
                     {\mbox{\boldmath$\scriptstyle{#1}$}}
                     {\mbox{\boldmath$\scriptscriptstyle{#1}$}}}}

\newcommand{\dom}{\mathsf{dom}}

\newcommand{\ip}[2]{{\left\langle {#1}, {#2} \right\rangle}}
\newcommand{\mip}[2]{\left\langle{#1}, {#2}\right\rangle_{\!\scriptscriptstyle{\text{M}}}}
\newcommand{\mat}[1]{\mathbf{{#1}}}
\newcommand\restr[2]{{
  \left.\kern-\nulldelimiterspace 
  {#1}\vphantom{\big|} \right|_{#2}}}

\newcommand{\R}{\mathbb{R}}

\newcommand{\B}{\mathcal{B}}
\newcommand{\C}{\mathcal{C}}
\newcommand{\D}{\mathcal{D}}

\newcommand{\F}{\mathcal{F}}

\newcommand{\J}{\mathcal{J}}

\newcommand{\bit}{\begin{itemize}}
\newcommand{\eit}{\end{itemize}}
\newcommand{\bdm} {\begin{displaymath}}
\newcommand{\edm} {\end{displaymath}}

\newcommand{\ProbMeasClass}{\mathcal{P}}
\newcommand{\ProbMeas}{\mathbbm{P}}
\newcommand{\SigmaAlg}{\mathcal F}

\newcommand{\GM}[2]{\mathcal{N}\!\left( {#1}, {#2}\right)}

\newcommand{\NoiseForcing}{\xi}

\newcommand{\obs}{\F}
\newcommand{\dobs}{{\vec{d}}}
\newcommand{\obsT}{{\vec{d}_{\scriptscriptstyle\text{obs}}}}
\newcommand{\ipar}{m}

\newcommand{\iparpr}{m_{\text{prior}}}

\newcommand{\ScoreFcn}{S}
\newcommand{\Nobs}{M}
\newcommand{\Expectation}{\mathbb{E}}

\newcommand{\like}{\pi_{\mbox{\tiny like}}}
\newcommand{\piobs}{\pi_{\mbox{\tiny obs}}}
\newcommand{\pinoise}{\pi_{\mbox{\tiny $\xi$}}}

\newcommand{\Cprior}{\mathcal{C}_{\text{\tiny{prior}}}}

\newcommand{\gbf}[1]{\boldsymbol{#1}}

\renewcommand{\vec}[1]{\gbf{#1}}

\usepackage{xspace}

\newcommand{\obsop}{\mathcal{B}}

\newcommand{\GN}{\ensuremath{\Gamma_{\!\!N}}}
\newcommand{\Exp}[1]{e^{{#1}}}
\newcommand{\GD}{\ensuremath{\Gamma_{\!\!D}}}
\newcommand{\Vg}{\mathscr{V}_{\!g}}
\newcommand{\V}{\mathscr{V}_{\!\scriptscriptstyle{0}}}

\newcommand{\LI}{\mathscr{L}}

\newcommand{\ipart}{m_{\scriptscriptstyle\text{true}}}
\newcommand{\grad}{\nabla}
\newcommand{\ut}[1]{\ensuremath{\tilde{#1}}}

\newcommand{\NEns}{{N_s}}

\newcommand{\CovFcn}{\rm{k}}
\newcommand{\dt}{{\Delta t}}

\newcommand{\priorcov}{\mat{\Gamma}_{\text{prior}} }
\newcommand{\priordensity}{\pi_{\text{prior}}}

\begin{document}

\maketitle

\vspace{-4.5cm}
\noindent{\color{red}\fbox{LLNL IM Release number:
    LLNL-JRNL-759502\qquad ANL Preprint \# ANL/MCS-P9140-1018}}
\vspace{4.5cm}

\begin{abstract}
This paper introduces a statistical treatment of inverse problems
constrained by models with stochastic terms. The solution of the
forward problem is given by a distribution represented numerically by
an ensemble of simulations. The goal is to formulate the inverse
problem, in particular the objective function, to find the closest
forward distribution (i.e., the output of the stochastic forward
problem) that best explains the distribution of the observations in a
certain metric.  We use proper scoring rules, a concept employed in
statistical forecast verification, namely energy, variogram, and
hybrid (i.e., combination of the two) scores.  We study the
performance of the proposed formulation in the context of two
applications: a coefficient field inversion for subsurface flow
governed by an elliptic partial differential equation (PDE) with a
stochastic source and a parameter inversion for power grid governed by
differential-algebraic equations (DAEs). In both cases we show that
the variogram and the hybrid scores show better parameter inversion
results than does the energy score, whereas the energy score leads to
better probabilistic predictions.
\end{abstract}

\section{Introduction\label{sec:concept}}

Inverse problems have been traditionally posed as inferring unknown or
uncertain parameters (e.g., coefficients, initial conditions, boundary
or domain source terms, geometry) that characterize an
underlying model from given (possibly noisy) observational or
experimental data~\cite{Vogel02,KaipioSomersalo05}. Such 
inverse problems governed by physics-based models, also referred to as data assimilation in the
meteorological and climate communities~\cite{Kalnay_B2003}, abound in
a wider range of application areas such as geophysics, cryosphere studies,
medical imaging, biochemistry, and control theory. Typically the
models governing these inverse problems are considered
deterministic. In reality, however, in addition to the inversion
parameter, these models involve other sources of uncertainties and randomness. For instance, the models have multiple uncertain coefficients or unknown or random source terms, parameters that are
not---or cannot be---inferred. Motivated by the need to account for
these additional uncertainties, researchers in recent years have shown a growing
interest in considering inverse problems governed by stochastic (or
uncertain) models, mostly in the
context of optimal control~\cite{Zabaras_2008,Borzi_2009,Rosseel_2012,Mattis_2015}. In this paper, we consider the inference
of parameters for stochastic models (described by differential
equations) and quantify the uncertainty associated with this inference.

\paragraph{Contributions}
This study introduces a methodology for the statistical treatment of inverse problems
constrained by physics-based models with stochastic terms. 
The salient idea of our approach is to express the  problem as finding the inversion
parameter for which the stochastic model generates a distribution that
best explains the distribution of the observations according to a loss-function we define. 
To this end, we formulate an  

inverse problem with a loss function given by suitable
statistical scoring 
metrics typically
employed in forecast verification
\cite{Thorarinsdottir10,Gneiting08,Gneiting05}. Proper scores allow for including in the inversion process a large range of statistical
features (e.g., spatial and/or temporal correlations and biases)
and, in this respect, are a significant departure from and improvement
over the traditional least-squares  misfit metrics. 
We also delve into the issue of how
different statistical scores affect the results of the inference
problem. This issue becomes important when we explore fitness functions
for multivariate distributions because one invariably needs to rely on
statistics that typically favor certain features over others, for example,
variance over correlations. 
Traditionally, the solution of the inverse problem is the parameter
field (or a distribution if working with statistical inverse problems)
that is the closest to the true parameter field in some metric, e.g.,
least-squares or statistical scores such as those proposed in this
work. However, one can also pose the problem as finding the parameter
field that generates the most accurate predictions in some statistical
sense. In other words, we are interested in how much more predictable
the model is after inference/inversion, not necessarily in the
goodness of fit. We also explore this alternative inversion paradigm
in this study and show that various proper scores or combination of
them can be used successfully in different inversion setups (i.e,
inverse problems governed by spatial differential equations and
time-dependent differential equations) to improve the model
predictability.

Another critical issue we address in this paper is the efficient computation of the numerical solution of the proposed statistical inverse problems. Namely, we propose a solution approach based on numerical optimization and provide the ingredients, in the form of gradient-based scalable optimization algorithms and adjoint sensitivities, that are needed to ensure the scalability of our methodology to large-scale problems. More specifically, 
to compute the most fit parameter field, we solve an optimization problem (implicitly) constrained by
the stochastic model with a quasi-Newton limited-memory secant
algorithm, BFGS updates for the inverse of the Hessian of the cost
function, and an Armijo line search. If the objective function is a
likelihood function and if prior information is available, then this
approach is equivalent to computing the maximum likelihood or maximum
a posteriori (MAP) 
point.

We derive adjoint-based expressions for efficient computation of the
gradient of the objective with respect to the inversion parameters. 
We illustrate our approach with two problems. The first is an inversion for the
coefficient field in an elliptic partial differential equation (PDE), interpreted as a subsurface flow
problem under a stochastic source field.
The second problem is represented by the inversion for a
parameter in a differential-algebraic system of equations (DAEs) under stochastic load terms. This model
can represent a power grid system with uncertain load behavior, which
induces small transients in the system, the goal here being to
determine the dynamic parameters within a quasi-stationary regime.

\paragraph{Problem formulation}
In what follows, let us consider that we have a mathematical model
expressed as $F(u,\ipar;\NoiseForcing)=0$, with states $u$ and
parameters $\ipar$ driven by a stochastic forcing $\NoiseForcing$  with
known probability law. Formally, such a model can consist of a standard stochastic
PDE on a domain $\D \subset \R^d$ ($d=1,2,3$) with suitable boundary
$\Gamma$; in this case $u$ is a function on $\D$ and $\NoiseForcing: \Omega \rightarrow \R^p$ ($p=2,3$) is defined by means of a probability space
$(\Omega,\SigmaAlg, \ProbMeas)$, where $\Omega$ is the sample space (the set of all possible events),
$\SigmaAlg$ is the $\sigma-$algebra of events, and
$\ProbMeas:\SigmaAlg\rightarrow[0,1]$ is a probability measure. We
assume that this probability space is completely defined.  We take 
  $\ipar$ to be a real-valued deterministic field, although extensions to
  random fields are also possible. Formally, the mathematical
  model can be defined as in~\cite{Hairer09} in the following form:
\begin{align}
  \label{eq:model:simulation}
  F(u(\cdot,\xi),\ipar(\cdot,\xi);\NoiseForcing(\cdot))=0 ~
  \textnormal{a.s.}\,,  ~\textnormal{with}~ F \in \D~\textnormal{a.e.},\,  u
  \in \overline{\D} \times \Omega \rightarrow \R,\, \ipar \in \R^d,\,
  \xi \in \R^{p}\,.  
\end{align}
We also assume that we have a suitable vector space 
with finite
stochastic moments when $F$ is a PDE. More details on the setup can be
found in~\cite{Gunzburger_2014a}. In this work $F$ is defined by a PDE
(see Sec. \ref{sec:poisson}) or an ordinary differential equation (ODE)/DAE (see
Sec. \ref{sec:num:power:grid}) system. We assume  the availability of
sparse observations $\obsT$ of the states $u$ corresponding to the
true parameters, which we denote by $\ipart$.

\bigskip

At the high level, the inverse problem formulation we use in this work
is standard: Given model~\eqref{eq:model:simulation} and observations
$\obsT = \F(\ipart) + \varepsilon_{\rm obs}$ such that
$F(u,\ipart;\NoiseForcing)=0$ a.s., find $\ipar$ that generates model
predictions that best explain the observations under a certain
metric. The function $\F(\ipart)$ is the so-called {\it
  parameter-to-observable map} whose evaluation involves the solution
of the given ODE/PDE followed by the application of an observation
operator $\obsop$, i.e., $\F(\ipart) = \obsop u(\ipar)$, where $u$
solves $F(u,\ipart;\NoiseForcing)=0$. The observations are subject to
known observational noise $\varepsilon_{\rm obs}$, which we assume to be a random
vector with known Borel probability measure $\piobs$, in addition to
and independent of the stochastic forcing $\xi$. A commonly used metric is the distance
between the observables predicted by the model $\F(\ipar)$ and the
actual observations $\obsT$. The metrics used in this work are
referred to as {\it statistical score functions} that compute the
fitness or a distance between the distribution of the observables
$\F(\ipar)$ and the set of validation data, namely, observations
$\obsT$. 
We denote these score functions by
$\ScoreFcn(\obs(\ipar),\obsT)$, where $\obs(\ipar)$ and $\obsT$
represent the model predictions and the observations, respectively.
We introduce and
discuss in detail such score functions in the next section.

\bigskip

Scores are positive functions that achieve their global minimum when
observations and model predictions are statistically indistinguishable.
For that reason scores have been used as loss or utility functions in order to
assess the level of confidence one has in the probabilistic model
prediction \cite{Kass_1995,Gneiting07}. Therein, the maximum score
estimation is introduced as a generalization of maximum likelihood
estimation. A likelihood function can be defined as
$\like(\obsT|\ipar) \propto
\exp\left(-\ScoreFcn(\obs(\ipar),\obsT)\right)$  as a
measure describing the relative plausibility of the parameter value \cite{Pfanzagl_1969}.
 Therefore, the inverse problem can be
formulated as finding  $\ipar^*$ such that
%
  \begin{align}
  \label{eq:model:optimization}
  \ipar^*=\arg
  \min_{\ipar}
\J(\ipar)
  ~\textnormal{ subject to }~F(u,\ipar;\NoiseForcing)=0\,\textnormal{ a.s. }
\end{align}
where $\NoiseForcing$ is known, and where, for instance,
$\J(\ipar)=-\log(\like(\obsT|\ipar))$.
In practice, we assume that we have access to a vector of $\Nobs$
observations $\obsT \in \mathbb{R}^\Nobs$ and can generate an
ensemble of $\NEns$ model predictions $\obs(\ipar)\in
\mathbb{R}^{\Nobs \cdot \NEns}$.

The objective in \eqref{eq:model:optimization} depends on
the observed data $\obsT$ via a single or multiple realizations, 
the numerical model
observables output $\obs(\ipar)$, and potentially explicitly on the parameter
$\ipar$. 
 In
this study we will follow the inverse Bayesian nomenclature with
$\ipar^*$ being the MAP point. We remark that the
optimization problem also depends on the parameter $\ipar$ implicitly
through the PDE or ODE/DAE constraint described by $F$.

\paragraph{Related work} Inverse problems with stochastic parameters are typically addressed
in a multilevel context. Here the stochasticity may come from reducing
the models and introducing a stochastic term that accounts for the
model reduction error \cite{Bal_2013,Lie_2017}. In other cases, the
additional stochastic/uncertain input is treated as a nuisance
parameter and an approximate premarginalization over this parameter is
carried out~\cite{NicholsonPetraKaipio18, KaipioSomersalo05,KaipioKolehmainen13}. In this study we consider
the case when stochasticity is inherent in the problem and we do not
have access to a deterministic version of a complete model.
In the optimal control
community, recent efforts have targeted moment matching between a stochastic
controlled PDE and observations,
\cite{Zabaras_2008,Borzi_2009,Rosseel_2012}. In most cases the loss
function is based on statistics of univariate or marginal
distributions. Various classical data discrepancy functions (utility/loss) for
inverse problems including Kullback-Leibler divergence are discussed \cite{Vogel_2002}.  Extending to multivariate settings is extremely
challenging because of the difficulty of accounting for complex
dependencies, the curse of dimensionality, and  the lack of 
order or rank statistics.

The scoring functions used in this study are
precisely addressing the multivariate aspect of the model and
observational data distributions.
A strategy similar to what we present here has been introduced in the
statistical community sometimes under the name of statistical
postprocessing or model output statistics. It consists of altering
the computational model probabilistic forecasts by postprocessing the ensemble
forecasts, and it tends to address the issue of bias and dispersion \cite{Gneiting05}.  
Most of these approaches are variations of Bayesian model averaging discussed in
\cite{Raftery05} and the nonhomogeneous regression  model
proposed in \cite{Gneiting05}. In these strategies the numerical model
or its parameters are not controlled; only its output is adjusted
\cite{Sloughter10,Lerch13,Baran14a,Baran14b,Scheuerer15,Thorarinsdottir10,Gel04}. In
the strategy proposed in this study, the model itself through its
parameters is part of the
control space. Therefore, our approach then can be interpreted as calibrating a generative
model~\cite{Bernton_2017} or model generator \cite{Semenov_1997,Maraun_2010}, where parameter $\ipar$ modulates the distribution of a
physical model simulator. In this context there are several strategies that aim
to minimize a certain distance between the generator and the
truth. The distance can be minimized between various statistics of the
generator outputs and the true data such as cumulative distribution
function, density, or moments. 
Those methods pertain to the class of minimum distance estimation in
which usual metric distances have been used such as the Chi-square,
the least-squares, or Kolmogorov-Smirnov. 
However, with increasingly complex models and complex distributions,
exact derivations of cost-functions becomes intractable and
approximation of distributions are obtained through strategies such as
Approximate Bayesian computations \cite{Marin_2012}. 
Many of these methods are emergent in the variational inference and machine learning communities \cite{zhang_2018},  where the Kullback-Leibler divergence tends to be the most popular combined with sampling algorithms. 
In the proposed work we propose use multivariate scoring metrics used
in statistical forecast evaluation field, which provide a
computational attractive, flexible, and
extensible alternative to existing metrics.

The remaining sections of this paper are organized as follows. We
begin by introducing the scoring functions and discuss the property
needed in order for the ansatz \eqref{eq:model:optimization} to be well posed in
Section~\ref{sec:scores}. In Sections \ref{sec:poisson}
and~\ref{sec:num:power:grid} we introduce a prototype elliptic
PDE-based model problem with application in subsurface flow and a
time-dependent problem driven by DAEs with application in power grid
modeling, respectively. We conclude in Section \ref{sec:conclusion}
with a discussion of the method, results, and its limitations.

\section{Scoring and metrics\label{sec:scoring}}
\label{sec:scores}
In statistics one way to quantitatively compare or rank probabilistic models is using
score functions.  A score function is a scalar metric $\ScoreFcn$ that
takes as inputs (i) verification data, in our inverse problem formulation the observations $\obsT$, and (ii) outputs from the
model to be evaluated, those outputs can be quantities describing the model (for instance parameters of probabilistic distribution) or model outputs, in the present work are the observables
subject to observational noise, namely, $\obs(\ipar)$, independent of
the stochastic forcing $\NoiseForcing$, and returns a scalar used to scoring or ranking verification outputs with respect
to the verification data. These scores are commonly
used in forecast evaluation where competing forecasts are
compared \cite{Thorarinsdottir_2013}. Scores are generally used in a negatively-oriented fashion, the smaller the score, the closer to the verification data is the model at stake.

While evaluating numerical model simulations, one aims to detect bias,
trends, outliers, or correlation misspecification.  
To create an objective function that can distinguish among different
modeling strategies, one needs appropriate mathematical
scoring metrics to rank them.  
Complex mathematical properties are required for consistent ranking of
the models \cite{Gneiting07}.

Score functions use different statistics to distinguish among different
(statistical) models. Moreover, in order to be able to distinguish among and consistently rank different models, score functions are required to have specific mathematical properties. \textit{Proper} score functions are widely used in statistics to ensure consistent ranking, for example in forecast evaluation, where competing forecasts are
compared~\cite{Thorarinsdottir_2013}. The following definition of proper scoring is from~\cite{Gneiting07}.  
\begin{definition}
  \label{def:proper}
  For the considered probability space $(\Omega,\SigmaAlg,\ProbMeas)$, a score $\ScoreFcn: \ProbMeasClass \times \Omega \rightarrow \mathbb{R}$ 
  is proper relatively to the class of probability measure $\ProbMeasClass$ iff
  \begin{align}
    \label{eq:proper:def}
    \Expectation_Y\{\ScoreFcn(P_Y, Y)\} \le \Expectation_Y\{\ScoreFcn( P , Y)\}\,, \forall P \in \ProbMeasClass\,,
  \end{align}
  where $\ProbMeasClass$ is a convex class of probability measures on $(\Omega, \SigmaAlg)$, and $P_Y$ is the probability distribution of $Y$.
\end{definition}
In other words, this definition states that a proper scoring rule prevents a score from favoring any
probabilistic distribution over the distribution of the verification
data. In addition, a proper score has two main desired features: (1) 
statistical consistency between the verification data and the model
outputs,  called {\it calibration}, and (2) reasonable dispersion in the model outputs, provided they are 
calibrated, which is referred to as {\it sharpness}.  In statistics,
the trend is to build proper scores that 
access simultaneous calibration and sharpness.

Various functions can be used to assess the error between the
verification data and  model outputs; however, scoring is not restricted to pointwise
comparison.  Methods to evaluate the quality of unidimensional outputs
are well understood~\cite{Gneiting_2014}; however, the evaluation of
multidimensional outputs or ensemble of outputs has been addressed in
the literature more recently \cite{Gneiting08,Pinson12,Scheuerer15b}
and remains challenging.  The {\it energy} and {\it variogram-based
  scores} (detailed below) are well suited to multiple
multidimensional realizations of a same model to be
verified. Therefore, in this paper we will focus on these scores and
combination of them. 

A widely accepted score is the continuous ranked probability
score (CRPS):
\begin{align}
  \label{eq:CRPS}
\ScoreFcn_\textnormal{CRPS}(P,y) = \int_{-\infty}^{\infty} \left(
         F_{X}(x) - \mathbbm{1}_{x<y} \right)^2 \mathrm{d} x\,,\quad \,
\end{align}
where $F_X$ is the cumulative distribution function (CDF) of $X\sim P$, $F_{X}(x) = \ProbMeas[X \le x]$, and
$\mathbbm{1}_{x<y}$ is the Heavyside function. 
The CRPS computes a distance between a full probability distribution
and a single deterministic observation, where both are represented by
their CDF. However, this score is only univariate and cannot be used
if the dimension of the observations is larger than one. 
In the following, we consider the energy score and the variogram-based score that are both scores expressed in a multi-dimensional context.

Moreover, closed forms of the scores are not always computable, consequently one uses Monte Carlo approximation of the scores by deriving them with samples from the predictive distribution of interest. For this reason, the energy and variogram-based scores will be computed using $\NEns$ samples in the following. 
 Approximated scores can then be expressed with discrete arguments as
  $\ScoreFcn:  \mathbb{R}^{\Nobs \cdot \NEns}  \times  \mathbb{R}^\Nobs  \rightarrow \mathbb{R} $, 
  which is applied to $\obs(\ipar)$, which is represented by $\NEns$ model prediction samples of dimension  $\Nobs$ and an observation or validation vector, $\obsT$, of size $\Nobs$.

\subsection{Energy score}
The energy score \cite{Gneiting08} is multivariate and proper. It 
generalizes CRPS \eqref{eq:CRPS} from univariate to multivariate and
can be expressed as 
\begin{align}
  \ScoreFcn(\dobs,\obsT) & =
  \label{eq:es:inf:dim:b}
  \Expectation_P\|\dobs^{a} - \obsT\| - \frac{1}{2}
  \Expectation_P\|\dobs^{a} - \dobs^{b}\| \,,~\obsT \sim
  P_{T}\,,~\forall \dobs^{a}, 
  \dobs^{b} \sim P\,,
\end{align}
where, in the context of this study, $\dobs=\obs(\ipar)$ which is considered a realization of probability distribution.
This score is sensitive to bias and variance discrepancy, but
potentially less sensitive to correlations; it will be denoted as the
\textit{ES-model}. 

In the probabilistic forecast context, scores can be used as a
loss function to fit probabilistic predictive distributions to
observations; for instance, see \cite{Scheuerer14}.  
Similarly to this idea, we propose to use statistical proper scores as objective functions in the underlying inverse problems. 
In this context, the score $\ScoreFcn$ could, for instance, be the
{\it energy score}; and if only samples from distribution $P$ are
available, it can be defined as follows:
\begin{align}
  \label{eq:energyscore}
  \ScoreFcn_{\mathrm{ES}}:=\ScoreFcn(\dobs,\obsT) = \frac{1}{\NEns}    \sum_{i=1}^{\NEns}
  || \dobs^{(i)} - \obsT|| 
  - \frac{1}{2\NEns^2}
  \sum_{i=1}^{\NEns}\sum_{j=1}^{\NEns}
  ||\dobs^{(i)}  -
  \dobs^{(j)}||, 
\end{align}
where $\NEns$ is the number of model prediction samples and
$\dobs^{(i)}=\{\obs(\ipar)\}^{(i)}=\B u^{(i)} (\ipar)$ are  model
predictions corresponding to $i$th sample of the stochastic model
forcing, $\NoiseForcing^{(i)}$, evaluated at parameter $\ipar$. Here
$\B$ 
is a linear observation operator that extracts measurements from
$u$. 

\subsection{Variogram score}
The variogram-based score \cite{Scheuerer15b} is multivariate, proper,
more sensitive to covariance (structure) but insensitive to bias. Its approximated sample version 
is given by
\begin{align}
  \label{eq:variogramscore}
  \ScoreFcn_{\mathrm{VS}}:=\ScoreFcn(\dobs,\obsT) = \sum_{i=1}^{\Nobs} \sum_{j=1}^{\Nobs} w_{ij}
  \left(|\obsT(i) - \obsT(j)|^p  -   \frac{1}{\NEns}
  \sum_{k=1}^{\NEns}|\dobs^{(k)}(i) -
  \dobs^{(k)}(j)|^p
  \right)^2\,,  ~ p > 0\,,
\end{align}
where we take $p=2$,  $\Nobs$ is the dimension of observations $\obsT$
(e.g., number of observational points in space), and $w_{ij}$ is a function of the distance 
of the position of observation $i$ and observation~$j$. In
other words we take differences between every observation and then of
the corresponding expectation of the scenarios.  If we
denote $\vec{\delta}_{ij} = \vec e_i-\vec e_j$, where $\vec e_i$ is the unit vector with
the $i^{th}$ component $1$, the preceding equation becomes
\begin{align}
  \label{eq:variogramscore2}
  \ScoreFcn(\dobs,\obsT) = \sum_{i=1}^{\Nobs} \sum_{j=1}^{\Nobs} w_{ij}
  \left(|\vec{\delta}^T_{ij} \; \obsT|^p  -   \frac{1}{\NEns}
  \sum_{k=1}^{\NEns}| \vec{\delta}^T_{ij} \; \B u^{(k)}|^p  \right)^2. 
\end{align}
This will be
referred to as the \textit{VS-model}.

\subsection{Discussion and other scores}
The energy score is known for failing to discriminate misspecified
correlation structures of the fields, but it successfully identifies fields
with expectation similar to the one of the verifying data.  
On the other hand, the variogram-based score fails to discriminate
fields with misspecified intensity, but it discriminates between
correlation structures \cite{Pinson13,Scheuerer15b}.  
Because of these different features and in order to discriminate fields
according to their intensity and correlation structure, we propose to use
a linear combination of the two scores, namely,  
\begin{align}
  \label{eq:combiscores}
  \ScoreFcn_{\mathrm{HS}}(\dobs,\obsT) = & \alpha \ScoreFcn_{\mathrm{ES}}(\dobs,\obsT) 
  + \beta   \ScoreFcn_{\mathrm{VS}}(\dobs,\obsT) \,,
\end{align}
where $\alpha>0$ and $\beta>0$ are problem specific. We will refer to
this hybrid score as the \textit{HS-model}. It is also a proper score
because any linear positive combination of proper scores remains a
proper score.

The score functions defined above are referred to as instantaneous scores
because they are functions of one verification data-point ($\obsT \in
\mathbb{R}^\Nobs$). If more than one 
verification sample is available---for example, if we have $n$ samples of
$\obsT^{(1,\dots,n)}=[\obsT^{(1)},\obsT^{(2)},\dots
  \obsT^{(n)}]^\top \in  \mathbb{R}^{n \times \Nobs}$ 
from the true distribution---then we can estimate the mean score
defined as follows:
\begin{align}
  \label{eq:mean:score}
  \ScoreFcn_n(\dobs,\obsT^{(1,\dots,n)})=\frac{1}{n}
  \sum_{i=1}^{n} \ScoreFcn(\dobs,\obsT^{(i)})\,.
\end{align}

In most cases, scores are used on verification data that are assumed
to be perfect. In practice, however,  observations are almost always
tainted with errors. Limited recent studies on forecast verification attempt to address this
issue \cite{Ferro17,Naveau18}. 
 Incorporating error and
uncertainty in the scoring setup is challenging. Analytical results
are tractable only in  particular cases such as linear or multiplicative noise with Gaussian and Gamma
distributions, respectively.  One way of tackling the observational
error is to assume some probability distributions for the observations
and the model outputs and to consider a new score defined as the
conditional expectation of the original score given the observations
\cite{Naveau18}.  Using the notations of Definition \ref{def:proper}, we can express
the corrected score  as  
$\ScoreFcn_{corr}(P,y)=\Expectation( \ScoreFcn(P,X)|Y=y) $, where $X$ represents
the hidden true state of the system.  
In practice,  to implement this method, one has to assume some
distribution for the $X$ and $Y$ and access an estimation of the
distribution parameters. If the errors are i.i.d., then their
contribution can be factored out when using the
score as a loss function. This is the case under consideration in this
study.

\paragraph{Statistical properties of scores} 
Approximated scores are asymptotically unbiased and consistent
(convergent in probability) by the virtue of the Law of Large Numbers.  
 Moreover, as discussed in the introduction, scores can be used as
loss functions, this procedure falls into the class of optimum contrast
estimation. Asymptotical results about their consistency of optimum contrast estimators can be found in
\cite{Pfanzagl_1969}. Strictly proper scoring rules as a 
 contrast function is discussed in \cite{Gneiting07}.   
 In the case of non-strictly propriety, which is also the case in our study, one
 may loose the uniqueness of the limit point.  
 Namely, under regularity assumptions, the optimum estimator would
 converge to a point that belongs to a set of optima of the proper
 score. Multivariate strictly proper scores for non-standard
 distributions are typically intractable, and hence, in this study we
 focus on proper scores, which are practical. Moreover, proper scores can be
 seen as divergence functions; however, they typically do not satisfy
 the triangular inequality. We will refer in text to distance or
 metric in this weaker sense.

%
\section{Model problems\label{sec:models}}
%
To probe the proposed statistical treatment of inverse problems
constrained by differential equations-based models with stochastic
inputs, we consider two model problems. The first is
a coefficient field inversion for subsurface flow governed by an
elliptic PDE with a stochastic input, in other words, a PDE-constrained model
problem (Section~\ref{sec:poisson}); the second is a parameter
identification problem for power grid governed by DAEs with stochastic
input, in other words, a DAE-constrained model problem
(Section~\ref{sec:num:power:grid}). For both problems, we generate
synthetic observations $\obsT$ by using one or more samples
$\NoiseForcing^{(i)} \sim \pinoise$, $i=1,\dots,\NEns$; where
$\pinoise$ is a known distribution. These samples then enter into the forward models with a
parameter considered the truth, $\ipart$. 
We then solve the optimization
problem~\eqref{eq:model:optimization} to obtain the maximum utility or
the maximum likelihood
 by evaluating
the likelihood function $\J(\ipar)=\ScoreFcn(\obs(\ipar),\obsT)$ and
for the MAP point by maximizing the a
posteriori probability density function $\J(\ipar)=\ScoreFcn(\obs(\ipar),\obsT)+ \mathcal{R}(\ipar)$
with the precomputed $\NEns$ scenarios $\NoiseForcing^{(i)}$ such that
$F(u^{(i)},\ipar;\NoiseForcing^{(i)}) = 0$, $i=1,\dots,\NEns$. To
solve the PDE-constrained model problem efficiently, in
Section~\ref{sec:derivatives} we derive the gradient of the objective
function $\nabla_\ipar\J(\ipar)$ with respect to
these fixed scenarios using adjoints. We remark that our calculations
 use classical Monte Carlo to estimate the solution of the
problem at hand; however, more sophisticated methods such as
higher-order \cite{Gunzburger_2017a} or multilevel \cite{Giles_2015}
Monte Carlo methods can be used to solve the underlying stochastic
PDE. 

\section{Model problem 1: Coefficient field inversion in an elliptic
  PDE with a random input} 
\label{sec:poisson}

In this section, we study the inference of the log coefficient field
in an elliptic partial differential equation with a random/stochastic
input. This example can model, for instance, the steady-state equivalent
for groundwater flows \cite{Mattis_2015}. For simplicity, we state the
equations using a deterministic right-hand side. We will then turn our
attention to the case where the volume source terms are
stochastic. To this end, consider the forward model
\begin{equation}\label{equ:poi}
  \begin{split}
    -\grad \cdot (\Exp{\ipar} \grad u) &= f \quad \text{ in }\D, \\
    u  &= g \quad \text{ on } \GD, \\
    \Exp{\ipar} \grad{u} \cdot \vec{n} &= h \quad \text{ on } \GN,
  \end{split}
\end{equation}
where $\D \subset \R^d$ ($d=2,3$) is an open bounded domain with
sufficiently smooth boundary $\Gamma = \GD \cup \GN$, $\GD \cap \GN =
\emptyset$.  Here, $u$ is the state variable; $f\in L^2(\D)$, $g\in
H^{1/2}(\GD)$, and $h\in L^2(\GN)$ are volume, Dirichlet, and Neumann
boundary source terms, respectively; and $\ipar$ is an uncertain
parameter field in $\mathcal{E} = \dom(\mathcal{A})$, where
$\mathcal{A}$ is a Laplacian-like operator, as defined in
~\cite{Stuart10,AlexanderianPetraStadlerEtAl16} and for
completeness repeated in Section~\ref{subsec:exp}. To state the weak
form of~\eqref{equ:poi}, we define the spaces,
\begin{align*}
    \Vg = \{ v \in H^1(\D) : \restr{v}{\GD} = g\}, \quad
    \V =  \{ v \in H^1(\D) : \restr{v}{\GD} = 0\},
\end{align*}
where $H^1(\D)$ is the Sobolev space of functions in $L^2(\D)$ with
square integrable derivatives.  Then, the weak form of \eqref{equ:poi}
is as follows: Find $u \in \Vg$ such that
\[
\ip{\Exp{\ipar} \grad{u}}{\grad{p}} = \ip{f}{p} + \ip{h}{p}_{\GN}, \quad
\forall p \in \V.
\]
Here $\ip{\cdot}{\cdot}$ and $\ip{\cdot}{\cdot}_{\GN}$ denote the
standard inner products in $L^2(\D)$ and $L^2(\GN)$, respectively.

In what follows we treat $f$ as a stochastic term, denoted by
$\xi$ for consistency, given by a two-dimensional heterogeneous
Gaussian process with known distribution $\pinoise$. We use the
instantaneous scores defined in Section~\ref{sec:scoring} for this model.

\subsection{Adjoint and gradient derivation}
\label{sec:derivatives}

We apply an adjoint-based approach to derive gradient information
with respect to the parameter field $\ipar$ for the optimization
problem~\eqref{eq:model:optimization} with
$\J(\ipar)=\ScoreFcn(\obs(\ipar),\obsT)+ \mathcal{R}(m)$, namely
\begin{equation}\label{eq:model:optimizationb}
  \min_{m \in \mathcal{E} } \ScoreFcn(\obs(\ipar),\obsT) +
  \mathcal{R}(\ipar)\,,
\end{equation}
where $\obs(m)$ corresponds to solving the forward
problem~\eqref{equ:poi} $\NEns$ times, and $\mathcal{R}(m)$, which will be
explicitly defined in the Computational experiment
Section~\ref{subsec:exp}, is a regularization/prior term.

The adjoint equations are derived through a Lagrangian
formalism~\cite{Troltzsch10}. To this end, the Lagrangian functional
can be written as
  \begin{equation}\label{eq:model:L}
    \LI(u,m,p):=  \ScoreFcn(\obs(\ipar),\obsT)  + \mathcal{R}(m) + \sum_{i=1}^{{\NEns}} \biggl [\ip{\Exp{m}\grad u^{(i)}}{\grad p^{(i)}}
  - \ip{\xi^{(i)}}{p^{(i)}} - \ip{p^{(i)}}{h}_{\GN}\biggr],
\end{equation}
where $p^{(i)} \in \V$ is the adjoint corresponding to state
$u^{(i)} \in \Vg$.  The formal Lagrangian formalism yields
that, at a minimizer of \eqref{eq:model:optimization}, variations of
the Lagrangian functional with respect to all variables vanish. Thus we have

\begin{subequations}
  \begin{align}
    \ip{\Exp{m} \grad u^{(i)}}{\grad \ut{p}} -
    \ip{\xi^{(i)}}{\ut{p}} - \ip{\ut{p}}{h}_{\GN} & = 0,     \label{eq:firststate}\\
    \ip{\Exp{m} \grad \ut u}{\grad p^{(i)}}
    + \ip{r^{(i)}}{\ut{u}} &= 0,         \label{eq:firstadj}\\
    \sum_{i=1}^{\NEns} \ip{\ut{m} \Exp{m}\grad u^{(i)}}{\grad p^{(i)}} &= 0,  \label{eq:firstcontrol}
  \end{align}
\end{subequations}
for all variations $(\ut{u}, \ut{m}, \ut{p}) \in \V \times \mathcal{E}
\times \V$ and $i = 1, \ldots, N_s$. Note that~\eqref{eq:firststate}
and \eqref{eq:firstadj} are the weak forms of the state and of the
adjoint equations, respectively.  The adjoint right-hand side
$r^{(i)}$ in strong form for the energy score~\eqref{eq:energyscore}
is
\begin{align} \label{eq:adjrhs}
  r^{(i)} = \frac{1}{2\NEns} \frac{\B^*(\B u^{(i)} - \obsT)}{|| \B u^{(i)} - \obsT||}
  - \frac{1}{\NEns^2} \sum_{j=1}^{\NEns}
  \frac{\B^*(\B u^{(i)} - \B u^{(j)})}{||\B u^{(i)}  - \B u^{(j)}||},
\end{align}
and for the variogram score~\eqref{eq:variogramscore} the $h^{th}$
component of  $r^{(i)}$ is
\begin{align} \label{eq:adjrhsvario}
  r_h^{(i)} = -\frac{4}{\NEns} \sum_{l=1}^M w_{lh} \C(u_h^{(i)},
  u_l^{(i)}) \; \B^* \vec{\delta_{lh}} \; \vec{\delta^T_{lh}} \; \B
  u^{(i)}, 
\end{align}
for $i = 1, \dots, \NEns$, and for $h = 1, \dots, M$. Here
\begin{alignat}{2}
  \C(u_h^{(i)}, u_l^{(i)}) &= |\obsT(h) - \obsT(l)|^2  -   \frac{1}{\NEns}
  \sum_{k=1}^{\NEns}|\underset{\B
    u^{(k)}(h)}{\underbrace{\{\obs(\ipar)\}^{(k)}(h)}} -
  \underset{\B u^{(k)}(l)}{\underbrace{\{\obs(\ipar)\}^{(k)}(l)}}|^2 \\&=
  |\vec{\delta}^T_{hl} \; \obsT|^2  -   \frac{1}{\NEns}
  \sum_{k=1}^{\NEns}| \vec{\delta}^T_{hl} \; \B u^{(k)}|^2, \nonumber
\end{alignat}
where $\B u^{(k)}(l)$ denotes the $l^{th}$ component of $\B u^{(k)}$,
namely, $\sum_{j=1}^M \B_{jl} u_j^{(k)}$.

The left-hand side in~\eqref{eq:firstcontrol} gives the
gradient for the cost functional~\eqref{eq:model:optimization}, which
is the Fr\'echet derivative of $\ScoreFcn(\obs(\ipar),\obsT)$ with
respect to $\ipar$. In strong form this is
\begin{align}
    \mathcal{G}(m) = \sum_{i=1}^{\NEns} \ip{\Exp{m}\grad
      u^{(i)}}{\grad p^{(i)}} + \mathcal{R}_m(m), \label{eq:grad_strong}
\end{align}
where $u^{(i)}$ and $p^{(i)}$ are solutions to the $i^{th}$ state and
adjoint equations, respectively, and $\mathcal{R}_m(m)$ is the
derivative of the regularization/prior term with respect to the
parameter $m$~\cite{Troltzsch10,BorziSchulz12}. The scaling of the
regularization term is problem specific and should be addressed
case-by-case.

We would like to add the following remarks about the adjoint problem:
(1) it is driven only by the derivative of the scoring functions with
respect to the forward solution; and (2) the forward and adjoint
problems share the same PDE operator, therefore the same solution
method can be applied to solve these PDEs. Computing the gradient
information via adjoints for large-scale PDE-constrained optimization
problems is imperative. Via an adjoint approach, the cost of the
gradient evaluation is one forward and one adjoint PDE solve per
optimization iteration~\cite{PetraStadler11}.

\subsection{Computational approach and cost}

  The inverse problems~\eqref{eq:model:optimization} are solved by using hIPPYlib (an
inverse problem Python library~\cite{VillaPetraGhattas18,VillaPetraGhattas16}).  It implements
state-of-the-art scalable adjoint-based algorithms for PDE-based
deterministic and Bayesian inverse problems. It builds on
FEniCS~\cite{DupontHoffmanJohnsonEtAl03, LoggMardalWells12} for the
discretization of the PDEs and on
PETSc~\cite{BalayBuschelmanGroppEtAl01a,BalayBuschelmanGroppEtAl09}
for scalable and efficient linear algebra operations and solvers needed for the solution of the PDEs.

The gradient computation technique presented in the preceding section
allows using state-of-the-art nonintrusive computational techniques of
nonlinear programming to solve the estimation
problems~\eqref{eq:model:optimization} efficiently for the energy and
variogram scores we propose, as well as any combination of them. More
specifically, we use a quasi-Newton limited-memory secant algorithm
with BFGS updates for the inverse of the
Hessian~\cite{Nocedal_book,BorziSchulz12} and an Armijo line search ~\cite{Nocedal_book} to solve~\eqref{eq:model:optimization} as an unconstrained optimization problem. 
This quasi-Newton solution approach is appealing since it can have fast local convergence properties similar to Newton-like methods without requiring Hessian evaluations and  it also converges from remote starting points as robust as a gradient-based algorithm.  In our computations the total number of quasi-Newton iterations was reasonably low, varying between $60$ and $160$. The implementation in hIPPYlib uses an efficient compact limited-memory representation~\cite{ByrdNocedalSchnabel_94_quasiNewtonRepres} of the inverse Hessian approximation that has reduced space and time computational complexities, namely, $O(|m|\times l)$, where $|m|$ denotes the cardinal of the discretization vector of $m$ and $l$ is the length of the quasi-Newton secant memory (usually taken as $O(10))$.

The computational cost per iteration is overwhelmingly incurred in the evaluation of the objective function in~\eqref{eq:model:optimization} and its gradient. For both the energy and variogram scores, the evaluation of the objective and its gradient requires $N_s$ forward PDE solves and adjoint PDE solves, respectively, to compute states $u^{(i)}$ in~\eqref{eq:firststate} and adjoint variables $p^{(i)}$ in~\eqref{eq:firstadj}. To achieve this, the projected states $d^{(i)}$ appearing in~\eqref{eq:energyscore} and~\eqref{eq:adjrhs} are stored in memory to avoid the expensive re-evaluations  of the PDEs and state projections for the computation of the score $S$ in~\eqref{eq:energyscore} and adjoint right-hand side in~\eqref{eq:adjrhs}. Similarly, for the variogram score, in the evaluation of the objective function we save the terms $ \C(u_h^{(i)}, u_l^{(i)})$ ($h,l=1,\ldots,M$) as a $M\times M$ matrix for each $i=1,\ldots,N_s$ and reuse them in the computation of the adjoint right-hand sides~\eqref{eq:adjrhsvario} during the objective gradient evaluation. This approach effectively avoids $N_s$ expensive re-evaluations  of the PDEs at the cost of $O(N_s M^2)$ extra storage. 

From \eqref{eq:energyscore} and~\eqref{eq:adjrhs} one can see that the computation of the energy score and its gradient also includes a $O(N_s^2 \cdot M)$ complexity term in addition to the forward and adjoint solves. A similar extra complexity term is present in the computation of the variogram score from~\eqref{eq:variogramscore} and its adjoint right-hand size from~\eqref{eq:adjrhsvario}.

Undoubtedly, the objective and gradient  computations  can be parallelized efficiently for both scores because of the presence of the summation operators. In particular, both scores allow a straightforward scenario-based decomposition that allows the PDE (forward and adjoint) solves to be done in parallel. Coupled with the (lower-level) parallelism achievable in hIPPYlib via DOLFIN and PETSc, this approach would result in an effective multilevel decomposition  with potential for massive parallelism and would allow tackling complex PDEs and a large number of scenarios. The quasi-Newton method based on secant updates used in this work can be parallelized efficiently, as one of the authors has shown recently~\cite{Petra_17_hiopdecomp}. 

We remark that a couple of potential nontrivial parallelization bottlenecks exist. For example, both the energy score and variogram score apparently require nontrivial interprocess communication in computing the right-hand side~\eqref{eq:adjrhs} of the adjoint systems  as well as in the computation of the double summation in the score itself. In this work we have used only serial calculations and deferred for future investigations efficient parallel computation techniques addressing such concerns.     

\subsection{Computational experiment} \label{subsec:exp}
In this section we present the numerical experiment setup for the forward
and inverse problems. 

{\it Forward problem:} For the forward problem~\eqref{equ:poi}, we
assume an unknown volume forcing, (i.e., $\xi \sim \pinoise$, with
known $\pinoise$) and
no-flow conditions on $\GN:=\{0,1\}\times (0,1)$, in other words, the
homogeneous Neumann conditions $\Exp{\ipar} \grad u \cdot \vec{n} = 0$
on $\GN$.  The flow is driven by a pressure difference between the
top and the bottom boundary; that is, we use $u = 1$ on $(0,1)\times
\{1\}$ and $u = 0$ on $(0,1)\times \{0\}$. This Dirichlet part of the
boundary is denoted by $\GD:=(0,1)\times\{0,1\}$.  In
Figure~\ref{fig:forwardprob}, we show the ``truth'' permeability used
in our numerical tests, and the corresponding pressure.

\def \pos {0.49\columnwidth}
\begin{figure}[t!]\centering
  \begin{tikzpicture}
    \node (1) at (0.0*\pos-0.77*\pos, 0.0*\pos){
    \includegraphics[height=.35\textwidth, trim=0 0 0 0, clip=true]{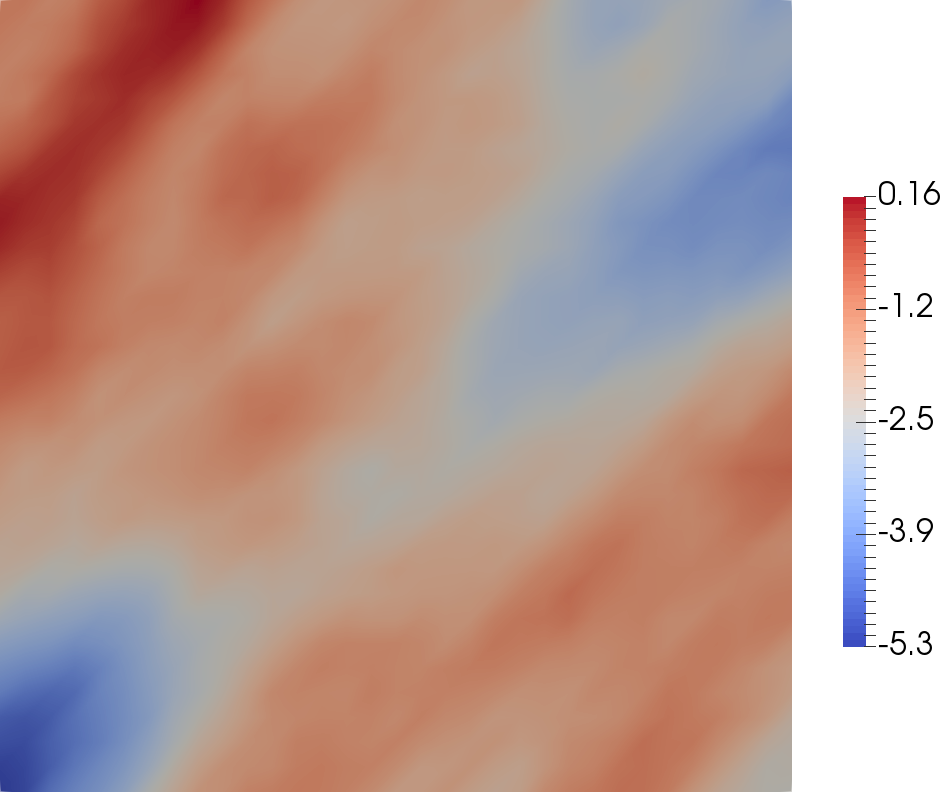}
    };
    \node (2) at (0.5*\pos-0.25*\pos, 0.0*\pos){
    \includegraphics[height=.35\textwidth, trim=0 0 0 0, clip=true]{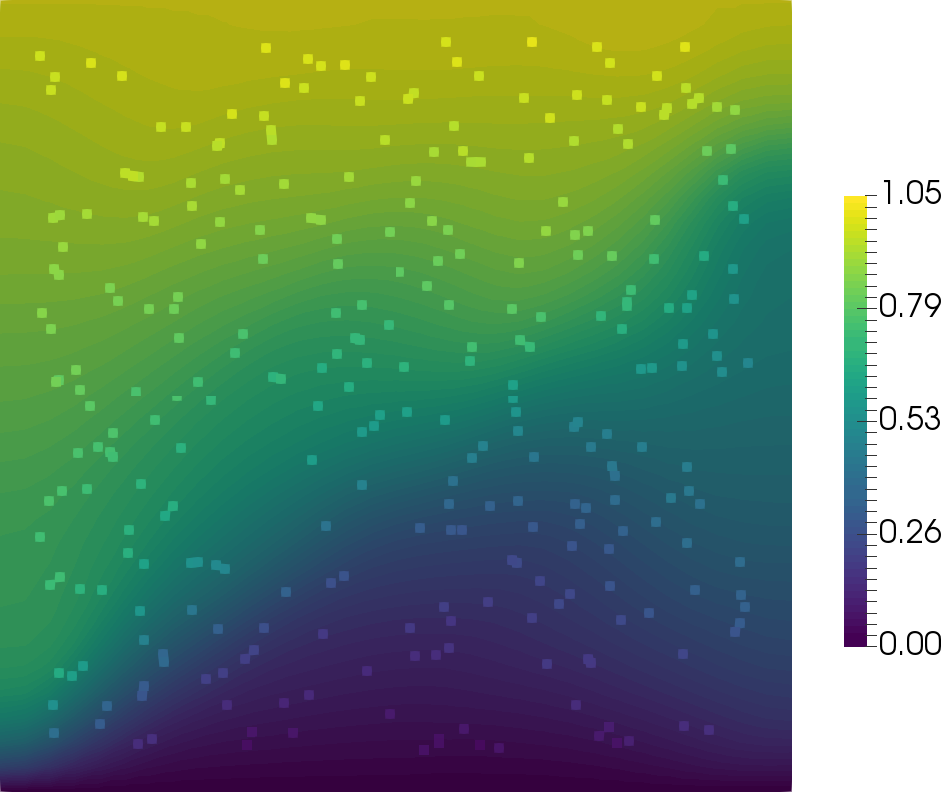}
    };
    \node at (0.0*\pos-1.15*\pos, -0.3*\pos) {\sf a)};
    \node at (0.7*\pos-0.83*\pos, -0.3*\pos) {\sf b)};
\end{tikzpicture}

\caption{Log permeability field $\ipart$ (a) and pressure field $u$
  obtained by solving the state equation with $\ipart$ (b). The dots
  show the location of observations $\obsT$.
}
\label{fig:forwardprob}
\end{figure}

\begin{figure}
\centering
\begin{tikzpicture}
   \node (1) at (0.0*\pos-0.85*\pos, 0.0*\pos){
    \includegraphics[width=.3\textwidth, trim=0 0 0 0, clip=true]{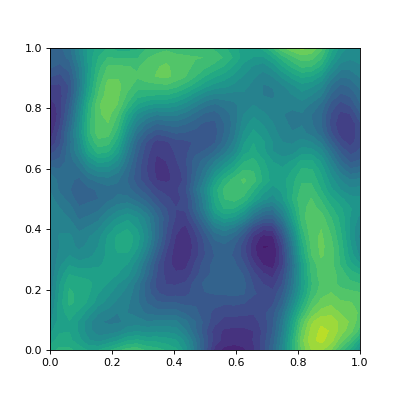}
    };
    \node (2) at (0.33*\pos-0.50*\pos, 0.0*\pos){
    \includegraphics[height=.3\textwidth, trim=0 0 0 0, clip=true]{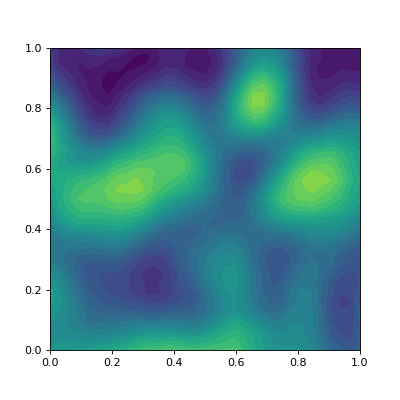}
    };
    \node (3) at (0.66*\pos-0.15*\pos, 0.0*\pos){
    \includegraphics[height=.3\textwidth, trim=0 0 0 0, clip=true]{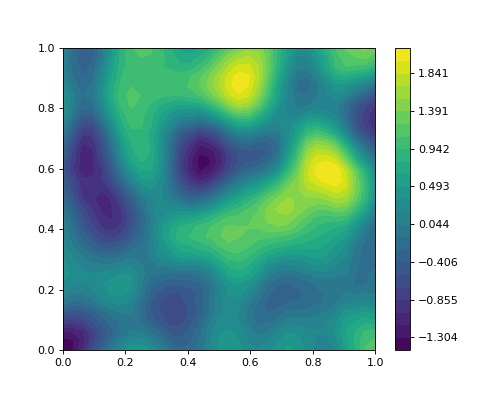}
    };
    \node (4) at (0.0*\pos-0.85*\pos, 0.0*\pos-0.7*\pos){
    \includegraphics[height=.3\textwidth, trim=0 0 0 0, clip=true]{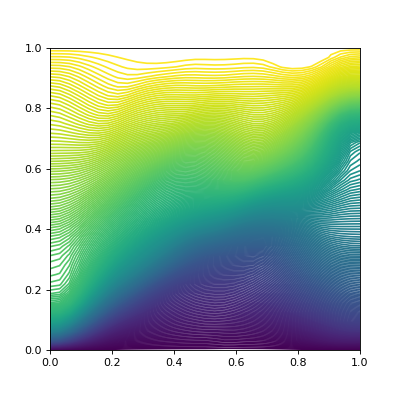}
    };
    \node (5) at (0.33*\pos-0.5*\pos, 0.0*\pos-0.7*\pos){
    \includegraphics[height=.3\textwidth, trim=0 0 0 0, clip=true]{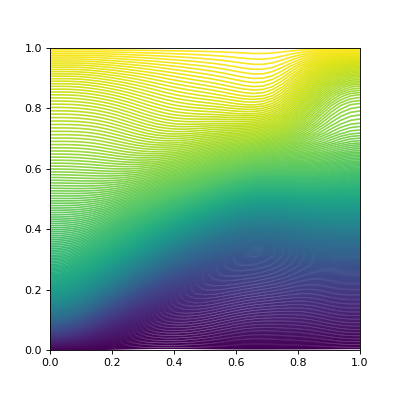}
    };
    \node (6) at (0.66*\pos-0.15*\pos, 0.0*\pos-0.7*\pos){
    \includegraphics[height=.3\textwidth, trim=0 0 0 0, clip=true]{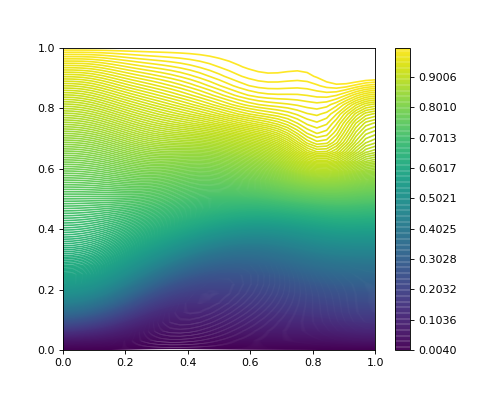}
    };
    \node at (0.0*\pos-1.0*\pos, -0.17*\pos) [white] {\sf $\xi$};
    \node at (0.0*\pos-0.3*\pos, -0.17*\pos) [white] {\sf $\xi^{(i)}$};
    \node at (0.0*\pos+0.3*\pos, -0.17*\pos) [white] {\sf $\xi^{(j)}$};
    \node at (0.0*\pos-0.95*\pos, -0.85*\pos) [white] {\sf $u(\ipart)$};
    \node at (0.0*\pos-0.3*\pos, -0.85*\pos) [white] {\sf $u^{(i)}$};
    \node at (0.0*\pos+0.3*\pos, -0.85*\pos) [white] {\sf $u^{(j)}$};
  \end{tikzpicture}

  \caption{
    \label{fig:subsurf:noise:output}
    Left column represents the pair of the noise realization $\xi$
    and solution $u(\ipart)$ used to generate the
    observations. 
  Two realizations
    of the forcing field ($\xi^{(i)},\xi^{(j)}$) (top) along with
    their corresponding solutions ($u^{(i)},u^{(j)}$) (bottom).}
\end{figure}
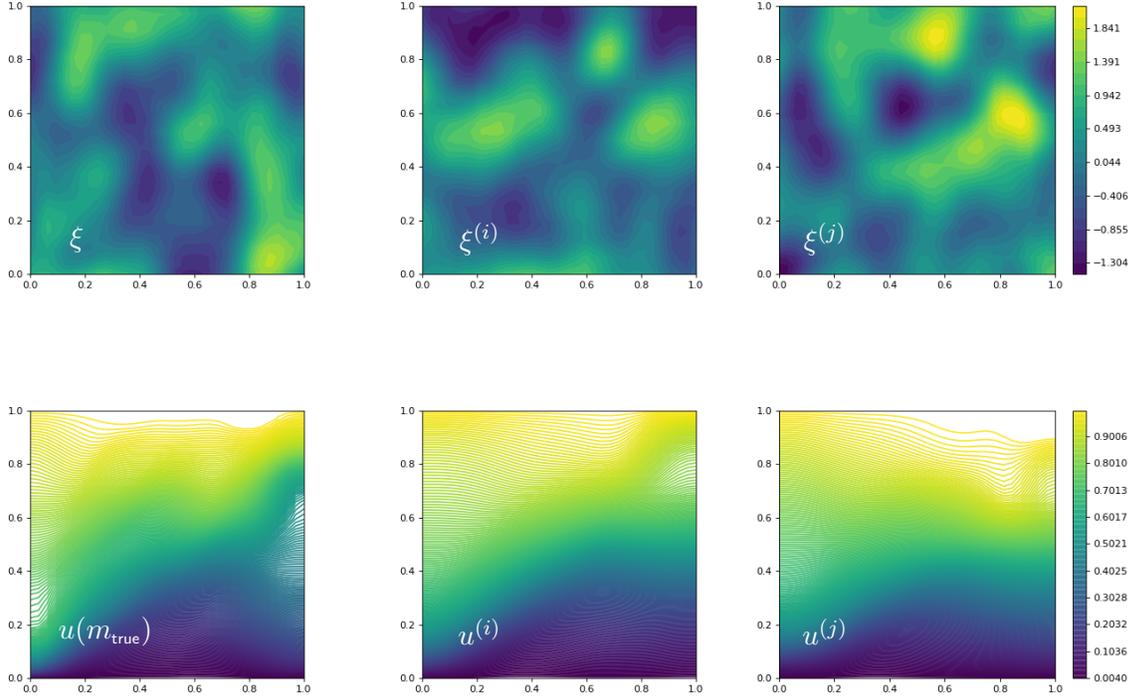

{\it The stochastic forcing term:} The stochastic volume forcing is
given by a two-dimensional heterogeneous Gaussian process with known
distribution $\pinoise$ defined by
\begin{align}
  \label{eq:elip:Noise}
  \xi(x,y) \sim {\cal N}\left(\mathbf{0},\CovFcn(h_x,h_y)\right)
,~ \CovFcn(h_x,h_y)=\sigma_\xi^2 \exp \left( -\frac{h_x^2}{\ell_{\Delta
        x}^2} -\frac{h_y^2}{\ell_{\Delta
        y}^2}  \right) + \delta_\xi I ,
  \end{align}
where $h_x=x-x'$, $h_y=y-y'$ and
$\mathrm{Cov}(\xi(x,y),\xi(x',y'))=\CovFcn(h_x,h_y)$. We choose
$\sigma_\xi=0.7$, $\ell_{\Delta x}=0.1875$, $\ell_{\Delta y}=0.1406$,
and $\delta_\xi=10^{-4}$. In Figure~\ref{fig:subsurf:noise:output} we
illustrate the forcing field and solution used to generate the
observations as well as two other realizations of the forcing field
along with their corresponding solutions.

{\it The observational noise:} We consider an uncorrelated (also
independent) observational noise with given distribution:
$\varepsilon_{\rm obs} \sim \mathcal{N}(0, \sigma^2 I)$, $\sigma^2 =
0.01$. As defined in the beginning, the observational noise is also
conditionally independent of forcing $\NoiseForcing$.

{\it The prior:} Following~\cite{Stuart10}, we choose the prior to
be Gaussian; that is, $\ipar \sim \GM{\iparpr}{\Cprior}$ is a prior
distribution, where $\iparpr$ is the mean and $\Cprior$ is the
covariance operator of the prior, modeled as the inverse of an
elliptic differential operator. To study the effect of the prior on
our results, we use an {\it informed} prior and the
{\it standard} prior, both built in
hIPPYlib~\cite{VillaPetraGhattas18}. The {\it informed} prior is
constructed by assuming that we can measure the log-permeability
coefficient at five points, namely, $N=5$, in $\D := [0,1]\times[0,1]$,
namely, ${\bf x}_1 = (0.1; 0.1)$, ${\bf x}_2 = (0.1; 0.9)$, ${\bf x}_3
= (0.5; 0.5)$, ${\bf x}_4 = (0.9; 0.1)$, ${\bf x}_5 = (0.9; 0.9)$, as
in~\cite{VillaPetraGhattas18}. This prior is built by using mollifier
functions
$$ \delta_i(x) = \exp\left( -\frac{\gamma^2}{\delta^2} \| x - x_i
\|^2_{\boldsymbol{\Theta}^{-1}}\right), \quad i = 1, \ldots, N.$$ The
mean for this prior is computed as a regularized least-squares fit of
the  point observations $x_i, i = 1, \ldots, N$, by solving
\begin{align}\label{equ:prior_mean_prob}
\iparpr = \operatornamewithlimits{\arg \min}_{m}
\frac{1}{2}\ip{m}{m}_{\widetilde{\mathcal{A}}} + \frac{p}{2}\ip{\ipart - m}{\ipart- m}_{\mathcal{M}},
\end{align}
where $\widetilde{\mathcal{A}}$ is a differential operator of the form
\begin{align}
  \widetilde{\mathcal{A}} = \gamma \nabla \cdot \left(
  \boldsymbol{\Theta} \grad \right) + \delta,
\end{align}
equipped with homogeneous natural boundary conditions,
$\mathcal{M} =
\sum_{i=1}^N \delta_i I$, and $\ipart$ is a realization of a
Gaussian random field with zero average and covariance matrix
$\mathcal{C} = \widetilde{\mathcal{A}}^{-2}$.  Above
$\boldsymbol{\Theta}$ is an s.p.d. anisotropic tensor, $\gamma$, and
$\delta > 0$ control the correlation length and the variance of the
prior operator; in our computations we used $\gamma = .1$ and $\delta
= .5$. The covariance for the informed prior is defined as $\Cprior
= \mathcal{A}^{-2}$, where $\mathcal{A} = \widetilde{\mathcal{A}} + p
\mathcal{M}$, with $p$ a penalization constant taken as 10 in our
computations. The {\it standard} prior distribution is
$\mathcal{N}(0, \Cprior)$, with $\Cprior =
\widetilde{\mathcal{A}}^{-2}$.

We note that the prior in finite dimensions is given by
\begin{align}
  \label{eq:prior_pdf}
  \priordensity(\ipar) \propto \exp\left[
    -\frac{1}{2}\mip{\ipar-\iparpr}{\priorcov^{-1}(\ipar-\iparpr)}\right],
\end{align}
where $\priorcov^{-1}$ is the discretization of the prior covariance
operator and $\mip{\cdot}{\cdot}$ is a mass weighted inner
product~\cite{Bui-ThanhGhattasMartinEtAl13,PetraMartinStadlerEtAl14}. In
the Bayesian formulation, the posterior is obtained as
$\pi(\ipar|\obs(\ipar),\obsT)\propto \like(\obs(\ipar),\obsT|\ipar)
\priordensity(\ipar)$. By taking the negative log of the posterior,
the objective in~\eqref{eq:model:optimization}  becomes 
\begin{align}
  \label{eq:obj:prior:bilaplacian}
  \J(\ipar))=\ScoreFcn(\obs(\ipar),\obsT) + \mathcal{R}(\ipar)
  \,.
\end{align}
where $\mathcal{R}(\ipar) =
\frac{1}{2}\mip{\ipar-\iparpr}{\priorcov^{-1}(\ipar-\iparpr)}$ and in
finite dimensional spaces $\mathcal{R}_\ipar(\ipar) = {\rm M}
\priorcov^{-1}(\ipar-\iparpr)$, where ${\rm M}$ is the mass matrix as
above. 

While the statistical
assumptions ease the computations, this objective can take a different
form under different functional likelihood or prior
expressions. However, the overall MAP finding procedure will broadly
follow the same steps.

\subsection{Results}

The aim of this computational study is two-pronged. On the one hand, our
goal is to invert for the unknown (or uncertain) parameter field with
the measure of success being the retrieval of a parameter field close
to the ``truth'': this is the traditional inverse problem
approach. On the other hand, we aim to generate accurate predictions in
some statistical sense; for example, we are interested in covering well the
multivariate distribution of the observables.  We will therefore carry out two
analyses: one  focused on the inverted parameter field  and one
on the model output (i.e., the observables). 
In each of the analyses, to assess the inversion quality of our approach, 
we will use 
the standard root mean square error (RMSE), the Structural
SIMilarity (SSIM) index \cite{Wang_2004}, and
the rank histogram as a qualitative tools.  The SSIM is the product of three terms (luminance,
contrast, and structure) evaluating respectively the matching of
intensity between the two datasets $a$ and $b$, the variability, and
the covariability of the two signals. In statistical terms,
luminance, contrast, and structure can be seen as evaluating the bias,
variance, and correlation between the two datasets, respectively. SSIM
is expressed as
$$SSIM(a,b) = \underbrace{\left(\frac{2\mu_{a}\mu_{b} + c_1}{\mu_{a}^{2} + \mu_{b}^{2} + c_{1}} \right)}_{\rm luminance} \underbrace{\left( \frac{2\sigma_{a}\sigma_{b} + c_2}{\sigma_{a}^{2} + \sigma_{b}^{2} + c_{2}}  \right)}_{\rm contrast} \underbrace{\left(\frac{\sigma_{ab} + c_3}{\sigma_{a}\sigma_{b} + c_{3}}   \right)}_{\rm structure},$$
where $\mu_{.}$, $\sigma_{.}$, and $\sigma_{..}$ respectively are the mean, standard deviation, and cross-covariance of each dataset, and $c_1$, $c_2$, and $c_3$ are constants derived from the datasets. 
The SSIM takes values between $-1$ and $1$. The closer to
$1$ the values are, the more similar  the two signals are in terms of intensity,
variability, and covariability. 
Researchers commonly also investigate the three components (luminance, contrast, variability) separately, as done hereafter. 

For a visual assessment of the statistical consistency between two
datasets in terms of probability
distributions~\cite{Anderson96,Hamill01} we use the rank histogram,
 an assessment tool often used in forecast verification.  This rank histogram gives us an
idea about the statistical consistency for the two datasets. The more
uniform the histogram is, the more statistically
consistent (i.e., sharp and calibrated) it is.

We solve the optimization
problem~\eqref{eq:model:optimization} with
$\ScoreFcn(\obs(\ipar),\obsT)$ as the energy and variogram scores and the
forward model $F(u,\ipar;\NoiseForcing)$ given by~\eqref{equ:poi}. To
understand the effect of the prior on the inversion results, for our
numerical studies, we consider two priors: an informed prior 
and a standard prior, 
as discussed in Section~\ref{subsec:exp}. In what follows, we discuss
the inversion results.

\paragraph{Comparison of MAP and true parameters}

\begin{figure}
  \centering
  \begin{tabular}{p{0.33\textwidth}p{0.33\textwidth}p{0.33\textwidth}}
    \\
    \centering
     $|\ipart - \ipar_\mathrm{guess}|$
      &\hspace{1.2cm}
      $|\ipart -
    \ipar_\mathrm{MAP}|$&\hspace{0.25cm}
     $|\ipart -
    \ipar_\mathrm{MAP}|$
    \\
    \includegraphics[width=0.41\textwidth]{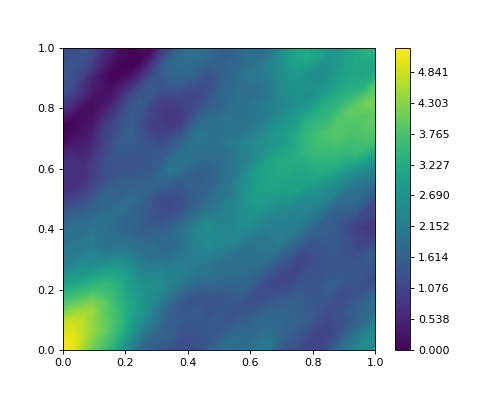}&\hspace{0.1cm}
    \includegraphics[width=0.33\textwidth]{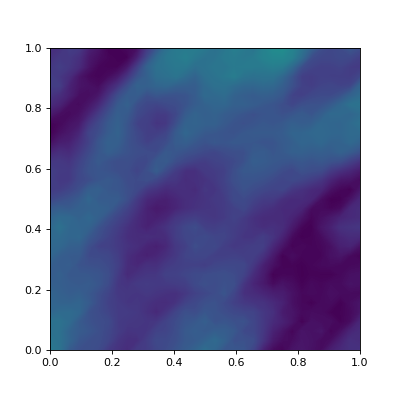}&\hspace{-0.75cm}
    \includegraphics[width=0.33\textwidth]{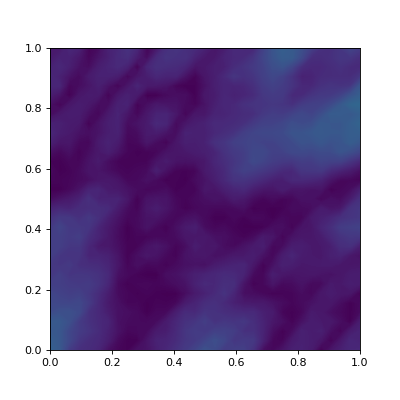}
    \\
    \centering
     initial guess &\hspace{1.2cm}
     ES-model (informed)&
     VS-model (informed)
    \\
    &\hspace{0.1cm}
    \includegraphics[width=0.33\textwidth]{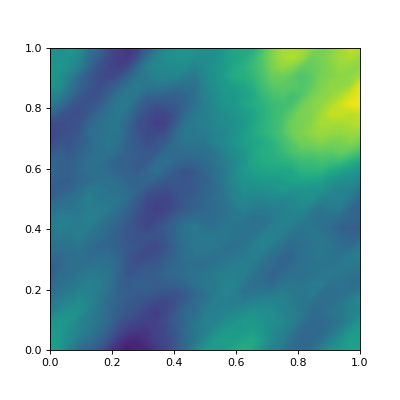}&\hspace{-0.75cm}
    \includegraphics[width=0.33\textwidth]{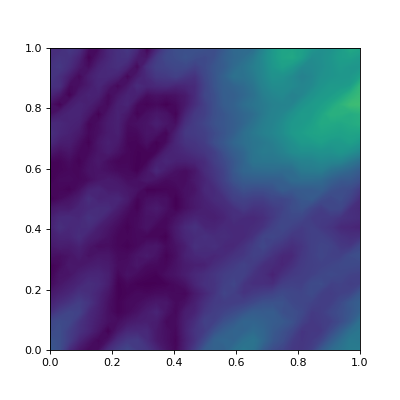}
    \\ 
    \centering &\hspace{1.2cm}
    ES-model (standard)&
    VS-model (standard)
  \end{tabular}
  \caption{Pointwise parameter field discrepancy $|\ipart -
    \ipar_\mathrm{MAP}|$ (left and center columns) and
    $|\ipart - \ipar_\mathrm{guess}|$ initial guess (top
    left) when using the informed (top) and standard (bottom) priors
    with the energy and variogram scores. The Monte Carlo sample size is
    64 in all panels. High discrepancy is indicated by light green,
    low discrepancy by dark blue.\label{fig:MAP:Discrepancy}}
\end{figure}

In Figure~\ref{fig:MAP:Discrepancy} we show the difference between the
true parameter $\ipart$ and the MAP estimate $m_{\rm MAP}$ for both the
informed (top) and standard priors (bottom), as well as for the
energy (left column) and variogram (center column) scores. In the
top-right panel we also show the initial guess for the optimization solver. The
results reveal that the VS-model objective 
   displays a stronger match between the MAP and
$\ipart$ than does the ES-model one. The HS-model (not shown in
   the figure) falls in
   between ES-model and VS-model as indicated in
   Table~\ref{tab:MAP:metrics}. The HS-model coefficients in
   \eqref{eq:combiscores} are 
   chosen to be $\alpha=0.1$ and $\beta=0.9$, with a better-informed
   choice possible but not fully explored in this study.
   Models with the informed prior exhibit smaller 
discrepancies than do the models with the standard prior. These results are
displayed in Table~\ref{tab:MAP:metrics}.
%


Table \ref{tab:MAP:metrics} shows the RMSE and SSIM and its 3
components, computed between the parameters $\ipart$ and $m_{\rm
  MAP}$.  As
expected, being in the informed prior case leads to smaller RMSE and better
SSIM for the inverted parameter, in other words, better overall
performance. 
Additionally, the contrast term, which is related to the
variance of each signal, is well captured by the three scores and by both types of priors (standard
and informed). In particular, we note that the use of standard priors
degrades the capture of the intensity of the parameters given by the
luminance term.  As expected, the VS-model and HS-model perform
better than the ES-model at capturing the covariance between the two
parameters $\ipart$ and $m_{\rm MAP}$ (structure term).
\begin{table}[!htbp]
\centering 
{\small \begin{tabular}{|p{1.2cm}|>{\raggedleft}p{1.3cm}|c|c|c|c|c|}
\hline
& {\bf  Samples} & {\bf  Luminance} & {\bf  Contrast} & {\bf  Structure} & {\bf  SSIM} & {\bf  RMSE}\\
\hline 
{\bf  Model} &  \multicolumn{6}{c|}{{\bf  (a) Informed prior}} \\
  \hline 
 \multirow{6}{=}{{\bf  ES} } & 1 & 0.847 & 1 & 0.698 & 0.591 & 1.137 \\
 & 4 & 0.801 & 0.989 & 0.797 & 0.631 & 1.134 \\
 & 8 & 0.824 & 0.995 & 0.777 & 0.637 & 1.108 \\
 & 32 & 0.803 & 0.996 & 0.765 & 0.612 & 1.155\\
& 64 & 0.793 & 0.992 & 0.757 & 0.596 & 1.176 \\
& 128 & 0.786 & 0.995 & 0.759 & 0.594 & 1.187 \\
\hline 
 \multirow{6}{=}{{\bf  VS} } & 1 & 0.825 & 0.997 & 0.755 & 0.621 & 1.124 \\
 & 4 & 0.955 & 0.976 & 0.868 & 0.809 & 0.696 \\
 & 8 & 0.966 & 0.967 & 0.859 & 0.803 & 0.66 \\
 & 32 & 0.98 & 0.948 & 0.848 & 0.788 & 0.617 \\
& 64 & 0.982 & 0.939 & 0.846 & 0.78 & 0.612 \\
& 128 & 0.987 & 0.947 & 0.855 & 0.799 & 0.574 \\
\hline 
 \multirow{6}{=}{{\bf  HS} } & 1 & 0.837 & 1 & 0.722 & 0.605 & 1.136 \\
 & 4 & 0.935 & 0.981 & 0.868 & 0.796 & 0.765 \\
 & 8 & 0.947 & 0.979 & 0.853 & 0.737 & 0.737 \\
 & 32 & 0.958 & 0.973 & 0.833 & 0.776 & 0.719 \\
& 64 & 0.960 & 0.966 & 0.83 & 0.769 & 0.716 \\
& 128 & 0.962 & 0.972 & 0.836 & 0.782 & 0.697 \\
\hline 
\hline 
{\bf  Model} &  \multicolumn{6}{c|}{{\bf  (b) Standard prior}} \\
\hline
 \multirow{6}{=}{{\bf  ES}} & 1 & -0.729 & 0.997 & 0.442 & -0.321 & 2.988 \\
 & 4 & -0.336 & 1 & 0.371 & -0.125 & 2.535 \\
 & 8 & -0.29 & 0.999 & 0.385 & -0.112 & 2.483 \\
 & 32 & -0.386 & 0.989 & 0.469 & -0.179 & 2.538 \\
& 64 & -0.283 & 0.985 & 0.469 & -0.131 & 2.43 \\
& 128 & -0.371 & 0.996 & 0.412 & -0.152 & 2.548 \\
\hline 
 \multirow{6}{=}{{\bf  VS} } & 1 & -0.828 & 0.991 & 0.527 & -0.432 & 3.146 \\
 & 4 & -0.547 & 0.999 & 0.447 & -0.245 & 2.737 \\
 & 8 & 0.268 & 0.998 & 0.363 & 0.097 & 2.017 \\
 & 32 & 0.412 & 0.995 & 0.46 & 0.189 & 1.8 \\
& 64 & 0.839 & 0.995 & 0.43 & 0.359 & 1.3 \\
& 128 & 0.863 & 0.999 & 0.402 & 0.346 & 1.292 \\
\hline 
 \multirow{6}{=}{{\bf  HS} } & 1 & -0.85 & 0.989 & 0.501 & -0.421 & 3.202 \\
 & 4 & -0.351 & 0.998 & 0.42 & -0.147 & 2.527 \\
 & 8 & -0.192 & 1 & 0.406 & -0.078 & 2.385 \\
 & 32 & -0.223 & 0.989 & 0.485 & -0.107 & 2.37 \\
& 64 & -0.023 & 0.992 & 0.465 & -0.011 & 2.194 \\
& 128 & -0.074 & 0.998 & 0.421 & -0.031 & 2.265 \\
 \hline 
\end{tabular} }
\caption{Quality of the reconstruction of the parameter field
  (i.e., the MAP point $m_{\rm MAP}$) measured by different metrics with
  informed (a) and standard (b) priors. The  {\bf
    Samples} column lists the number of Monte Carlo samples used to
  approximate the stochastic right-hand side was. The  {
    \bf Luminance} column shows the consistency in terms of intensity of the two signals;  
    the 
  {\bf Contrast} column represents the matching of variance of the two signals; and the {\bf Structure} column
  shows the covariance matching between the two signals. The {\bf SSIM} column---the product of the luminance, contrast, and structure---is a global measure of consistency of the two studied quantities. 
One expects the SSIM and its factor components to be as close to 1 as possible.
  The last column, {\bf RMSE}, shows the root mean squared errors between the two signals, one expects the RMSE to be as small as possible.
  This table shows that in the proposed setup of informed priors enables better results in terms of SSIM and RMSE and that the VS-model tends to provide a better matching between the true and estimated parameter field. The variance of each signal tends to be well captured by all models. }
  \label{tab:MAP:metrics}
\end{table}
%

\paragraph{Comparison between $\obs(\ipar)$ and $\obsT$}

\begin{figure}
  \centering
  \begin{tabular}{ccc}
\includegraphics[scale=.2]{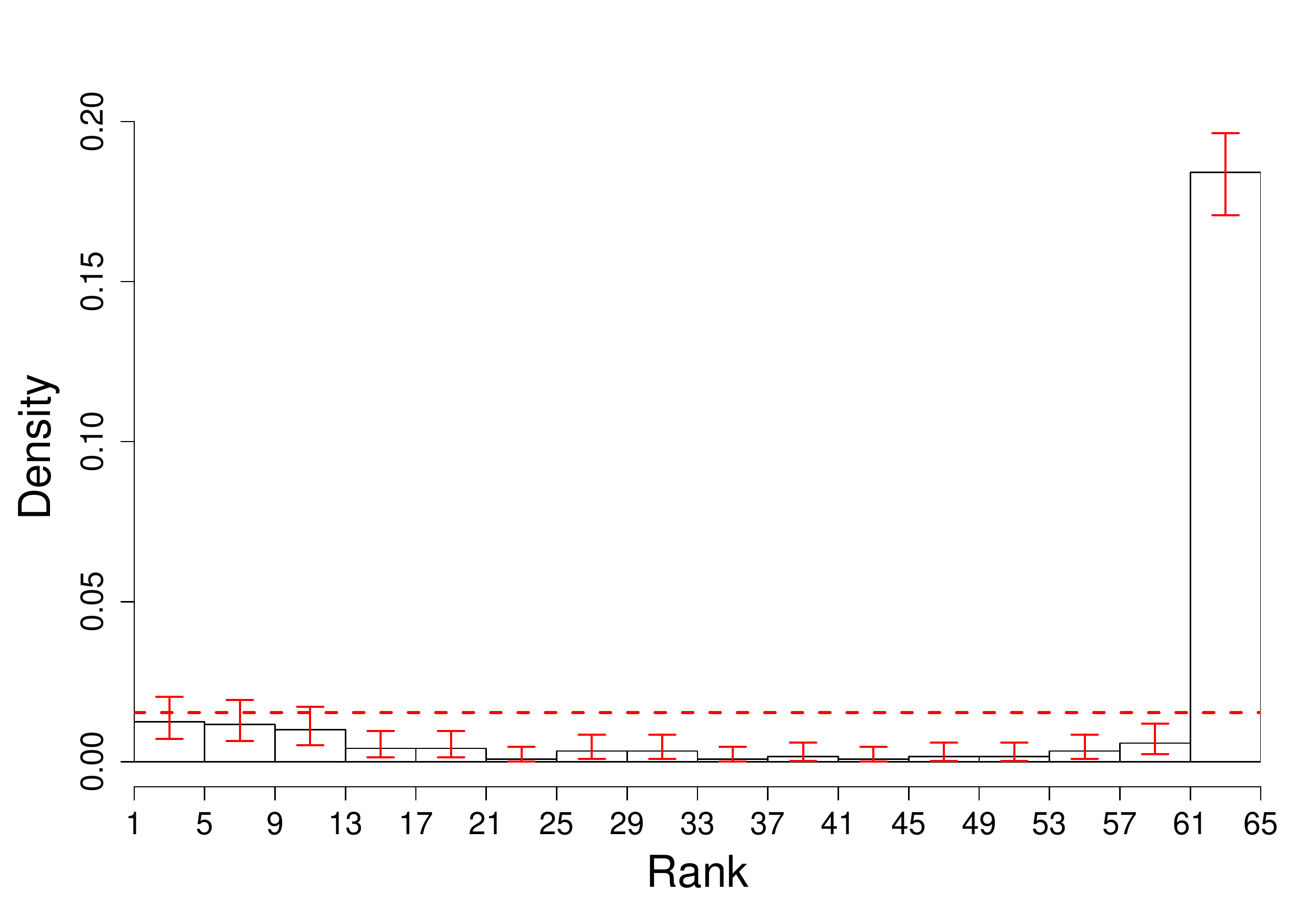}&
\hglue-.2cm
\includegraphics[scale=.2]{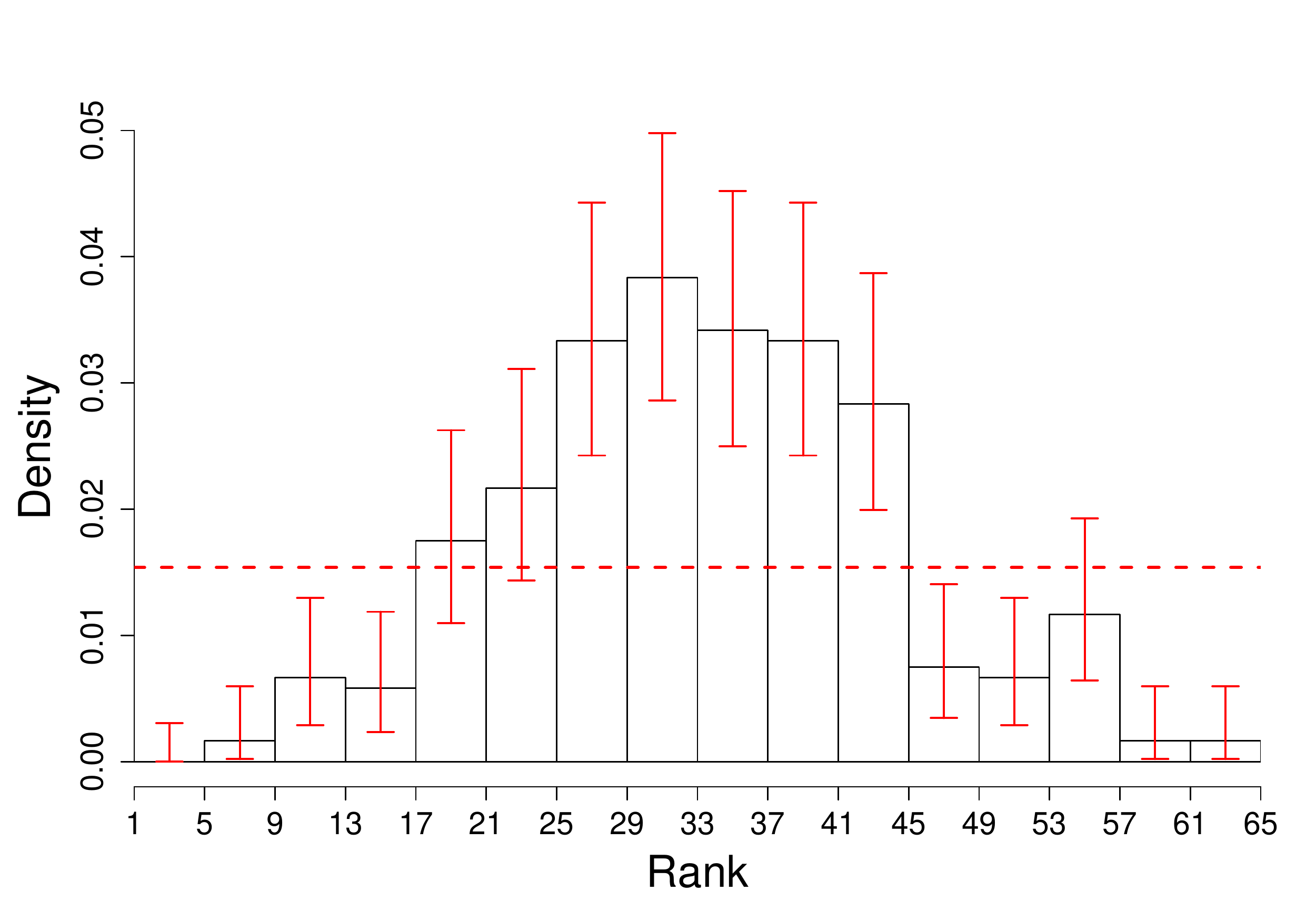}&
\hglue-.2cm
\includegraphics[scale=.2]{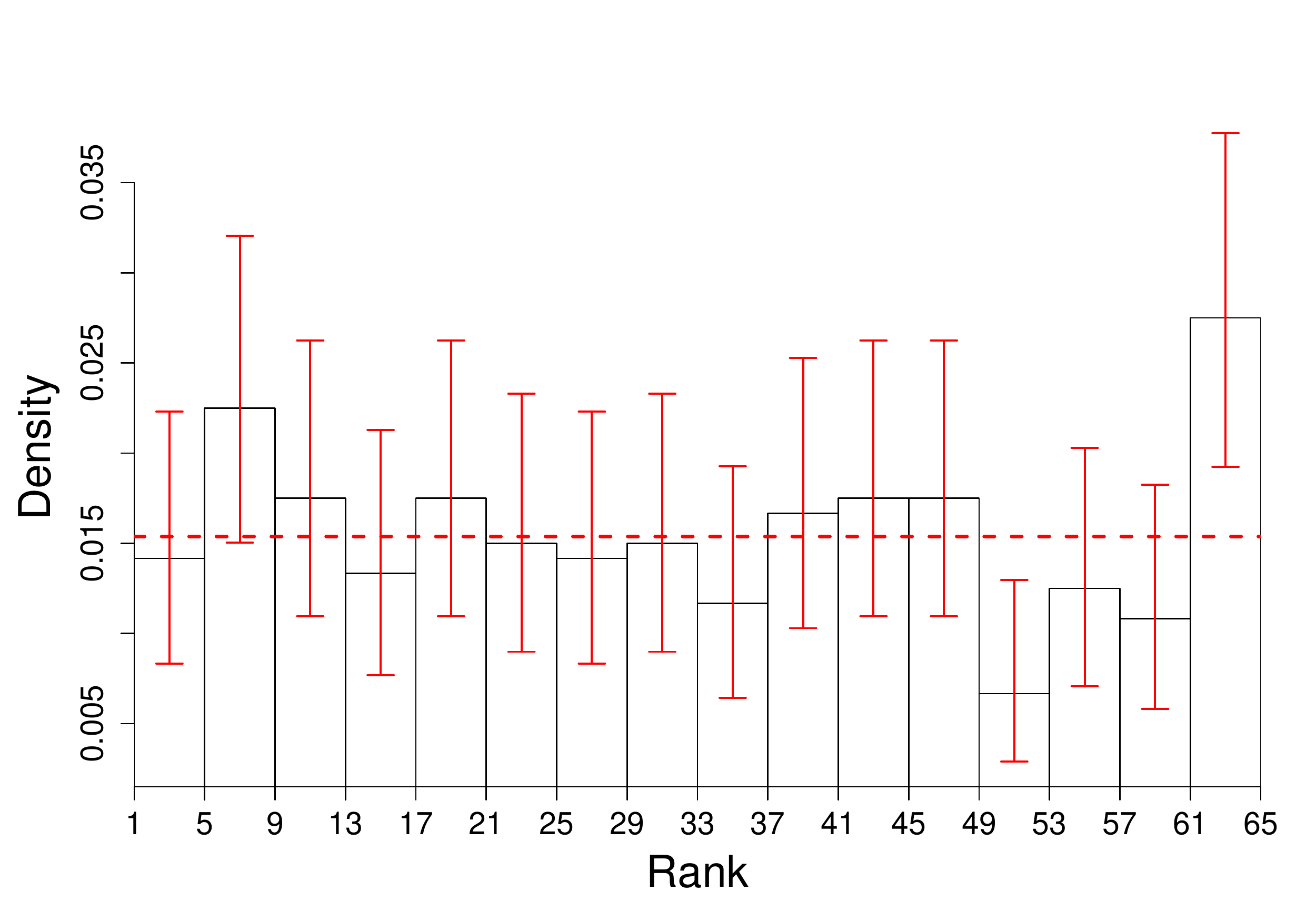} \\
initial guess
&
ES (informed prior)
&
ES (standard prior)\\
&
\hglue-.2cm
\includegraphics[scale=.2]{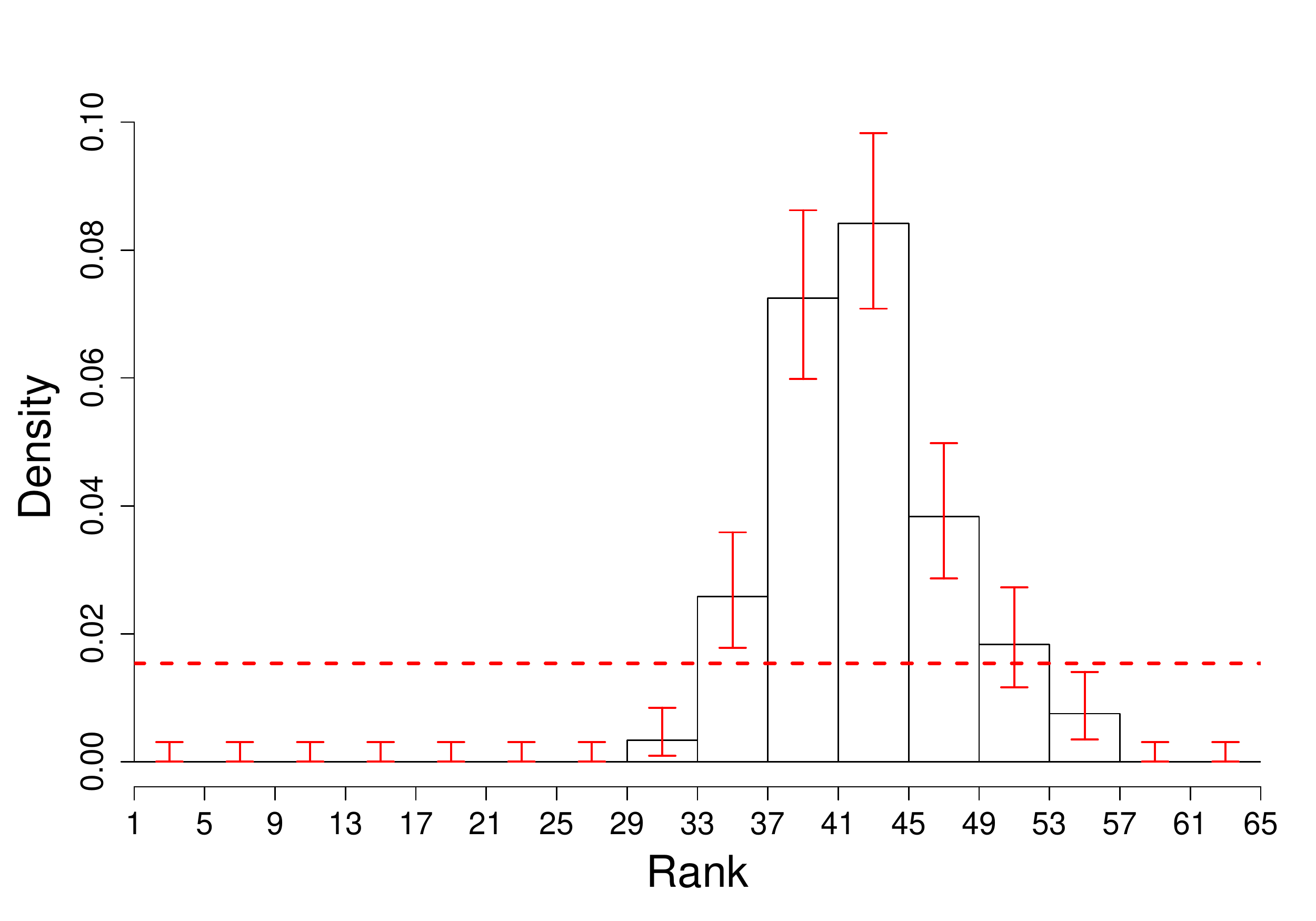}&
\hglue-.2cm
\includegraphics[scale=.2]{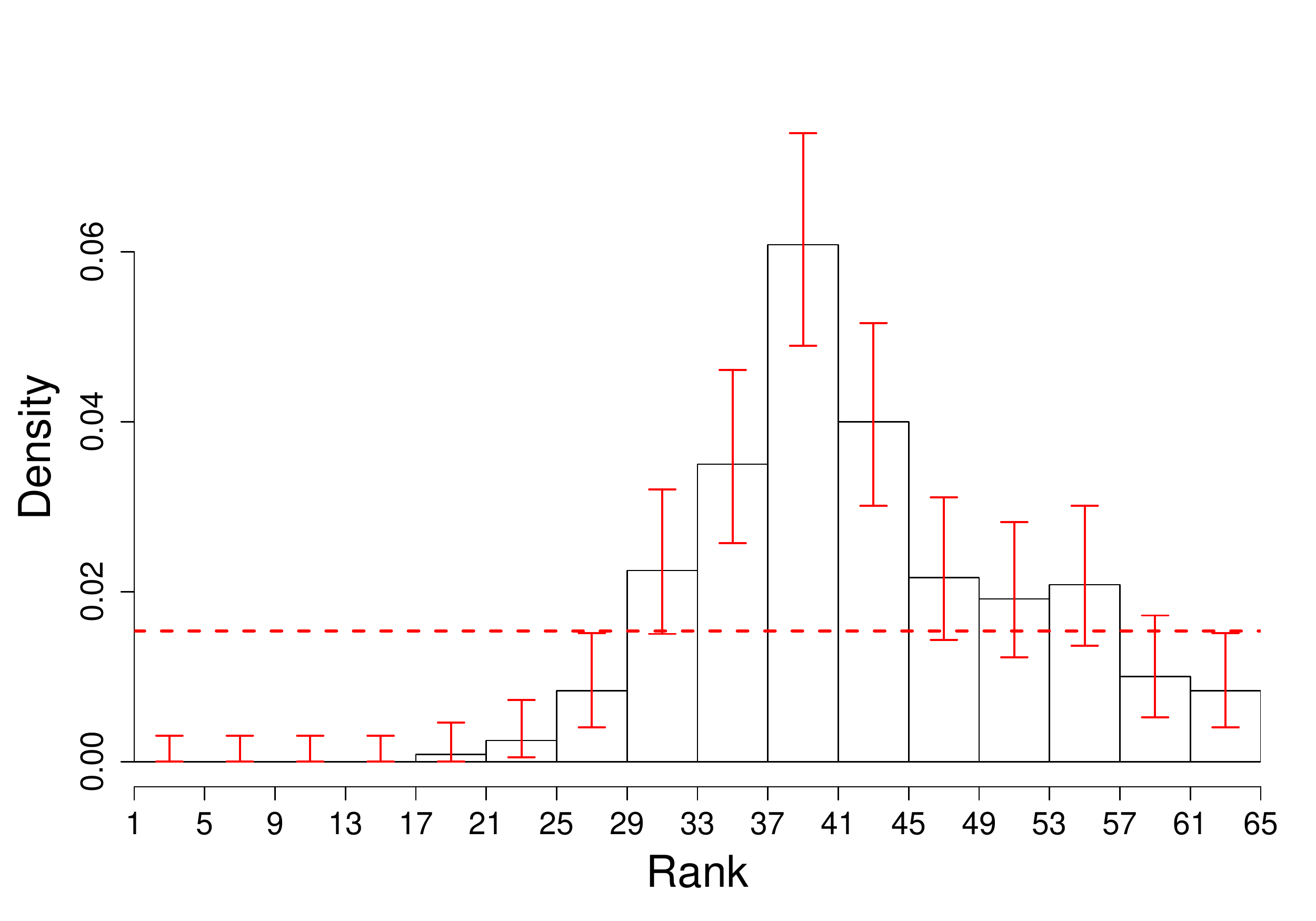} \\
&
VS (informed prior)
&
VS (standard prior)\\
&
\hglue-.2cm
\includegraphics[scale=.2]{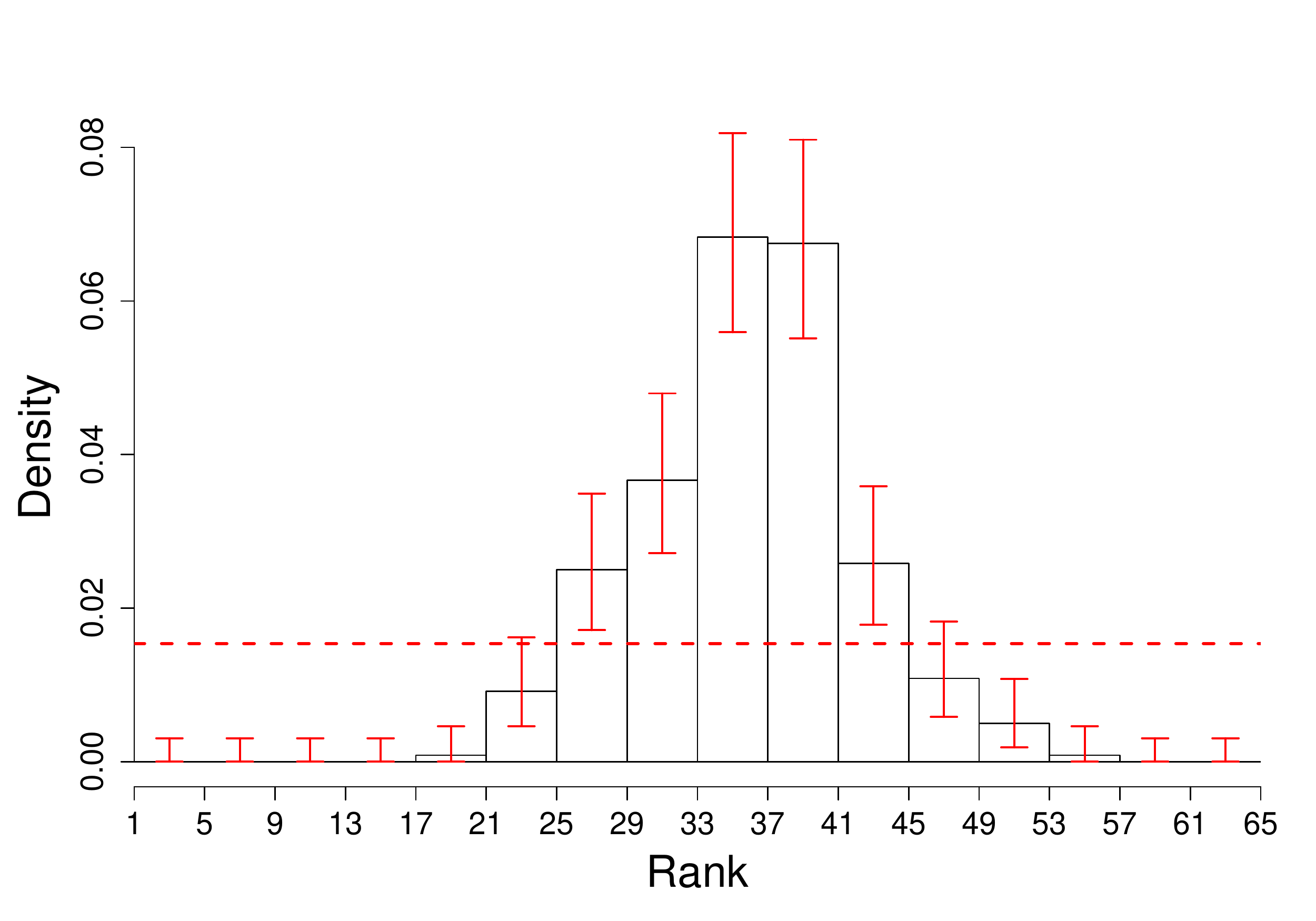}&
\hglue-.2cm
\includegraphics[scale=.2]{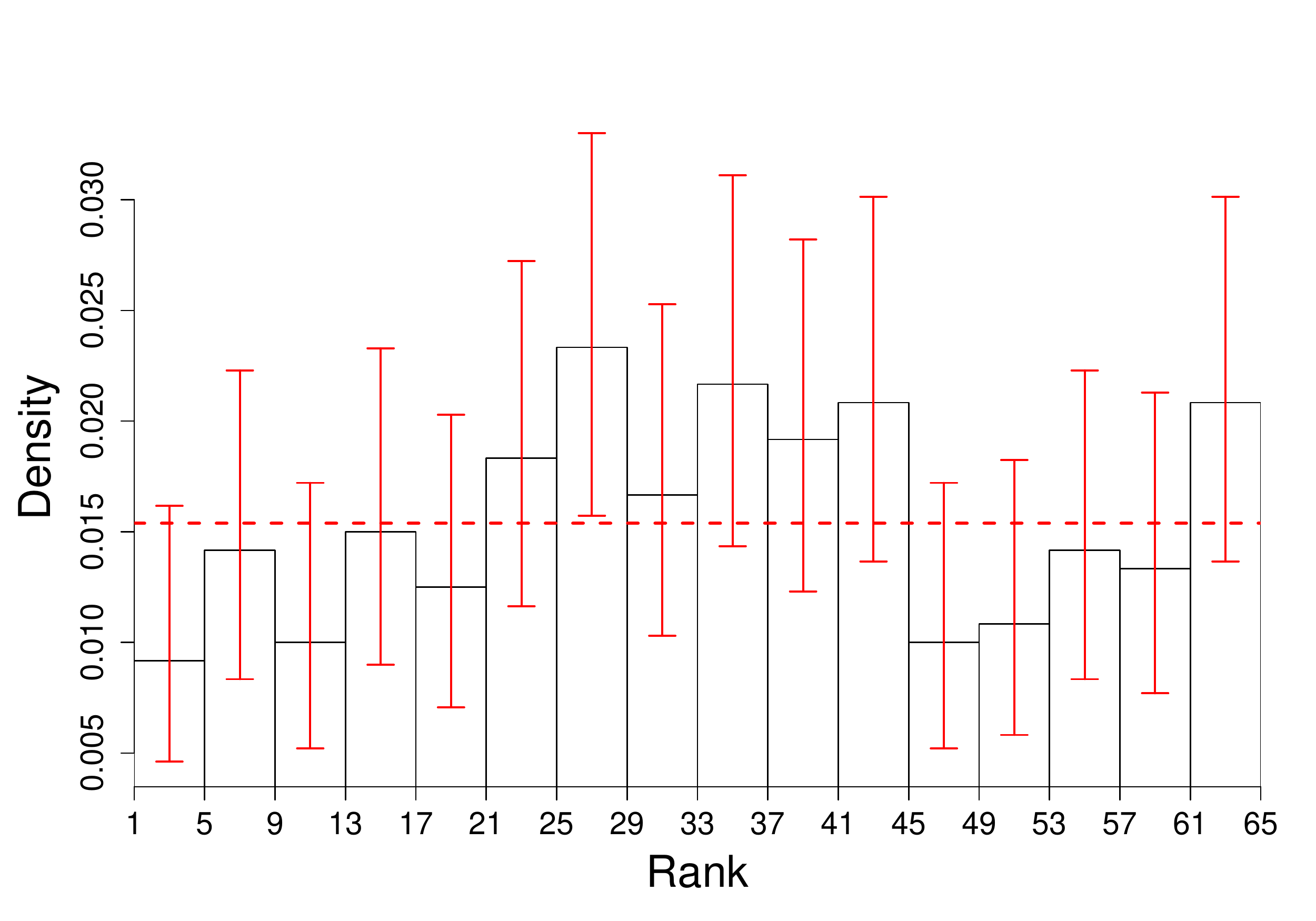}\\
&
HS (informed prior)
&
HS (standard prior)\\
\end{tabular}
\caption{Rank histogram between $\obsT$ and $\obs(\ipar)$. From left
  to right: guess, $\obs(\ipar)$ obtained by using the informed prior,
  and $\obs(\ipar)$ obtained by using the standard prior. From top to
  bottom: ES-model; VS-model; HS-model. The horizontal line is the
  perfect uniform histogram that represents a perfect match between
  $\obsT$ and $\obs(\ipar)$. Red whiskers show $\%95$-confidence
  intervals associated with the estimated count histogram. The closer
  to the uniform histogram, the better the consistency between $\obsT$
  and $\obs(\ipar)$. Histograms are obtained for simulations with 64
  samples.}\label{fig:obs:PIThisto}
\end{figure}

To assess the quality and statistical properties of the observables
generated by the model, in Figure~\ref{fig:obs:PIThisto} we show the
rank histograms reflecting the statistical consistency between the
true observables $\obsT$ and the \textit{generated} ones $\obs(\ipar)$. The
results show that the standard priors (right row) provide a better
calibration between $\obsT$ and $\obs(\ipar)$ than do the informed priors
(central row).  Additionally, the results show that the ES-model 
(top row) generates calibrated $\obs(\ipar)$. This is not unexpected since
the energy score is known for discriminating between the intensity of the
signals it compares~\cite{Pinson13}. The VS-model
(center row) does not present good calibration results. This result
is not unexpected either since the variogram score is known for not
capturing the intensity of the signals it
compares~\cite{Scheuerer15b}. We remark that the observables from the
VS-model show some overdispersion (bell-shaped histogram).  The
hybrid score seems
to take advantage of the properties of the energy score in terms of
calibration and thus appears to be a good compromise between the
ES-model and the VS-model. We note that the rank histograms do not
assess the correlation structure of the data. In that context we
investigate in the following indexes that measure the spatial data
structure.

\begin{figure}
\centering
  \begin{tabular}{ccc}
   \multicolumn{3}{c}{RMSE} \\[2pt]
\includegraphics[scale=.15]{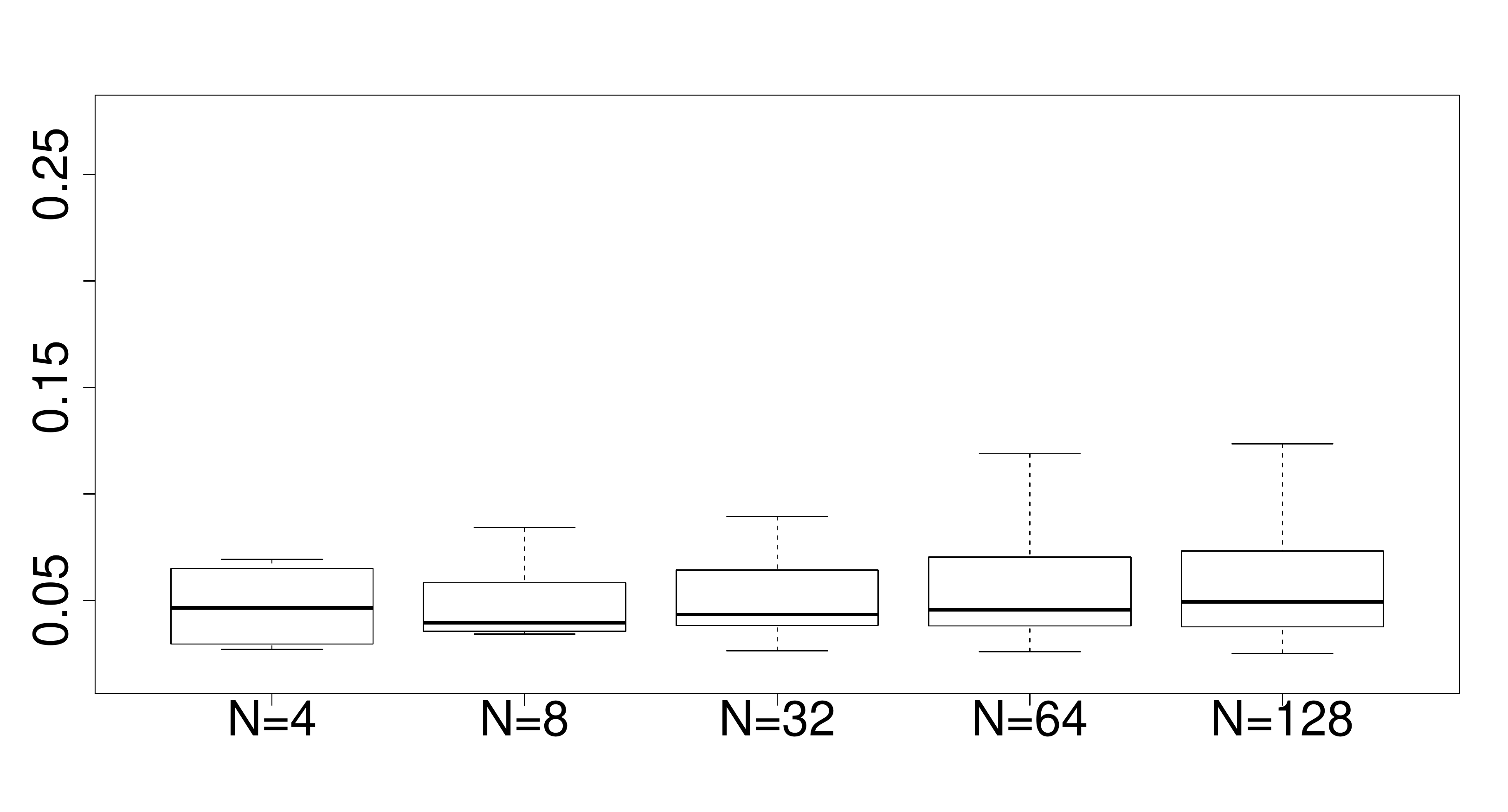}
& \includegraphics[scale=.15]{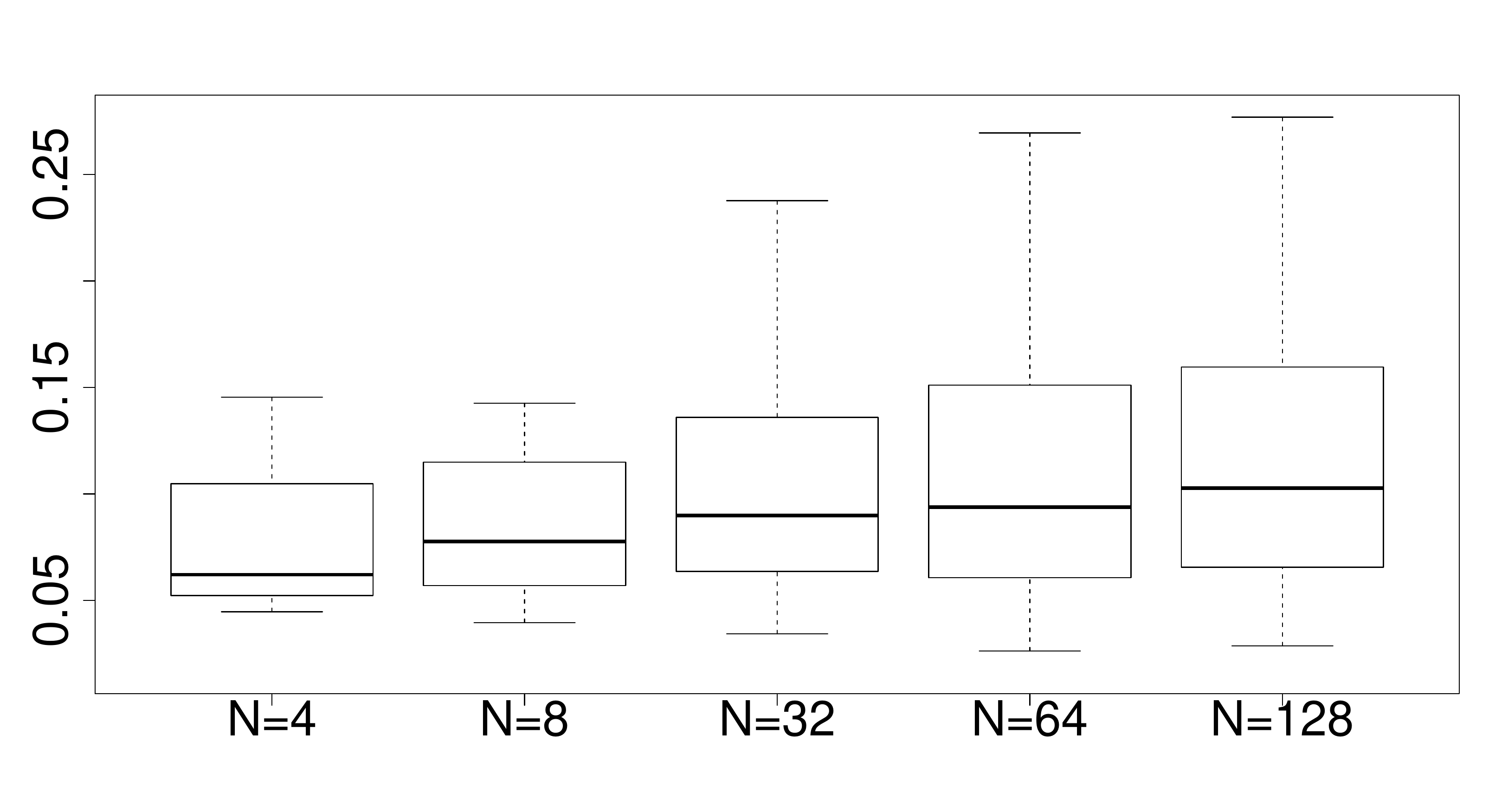}
& \includegraphics[scale=.15]{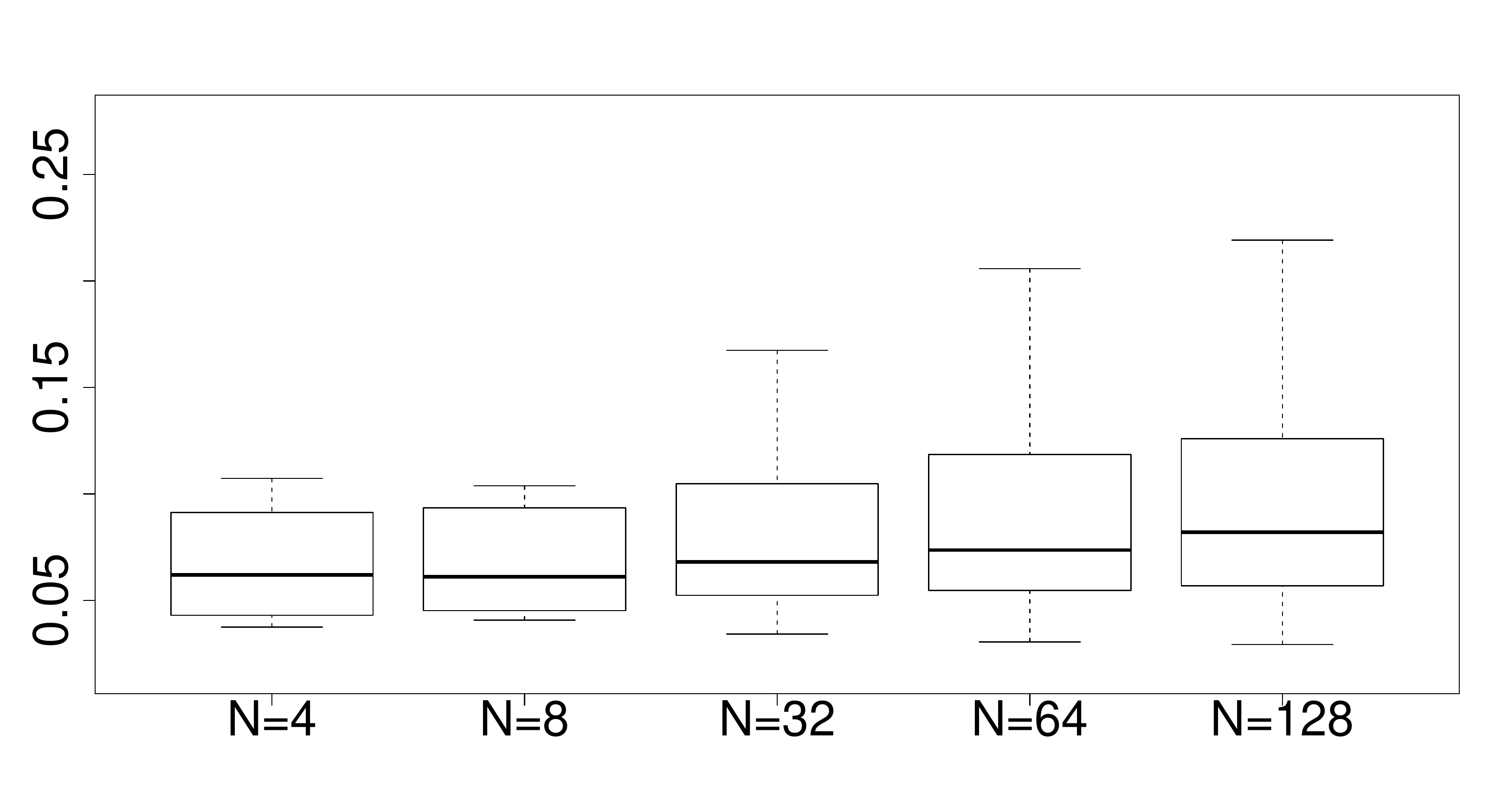} \\
  ES (standard prior) &  VS (standard prior) &  HS (standard prior) \\
\includegraphics[scale=.15]{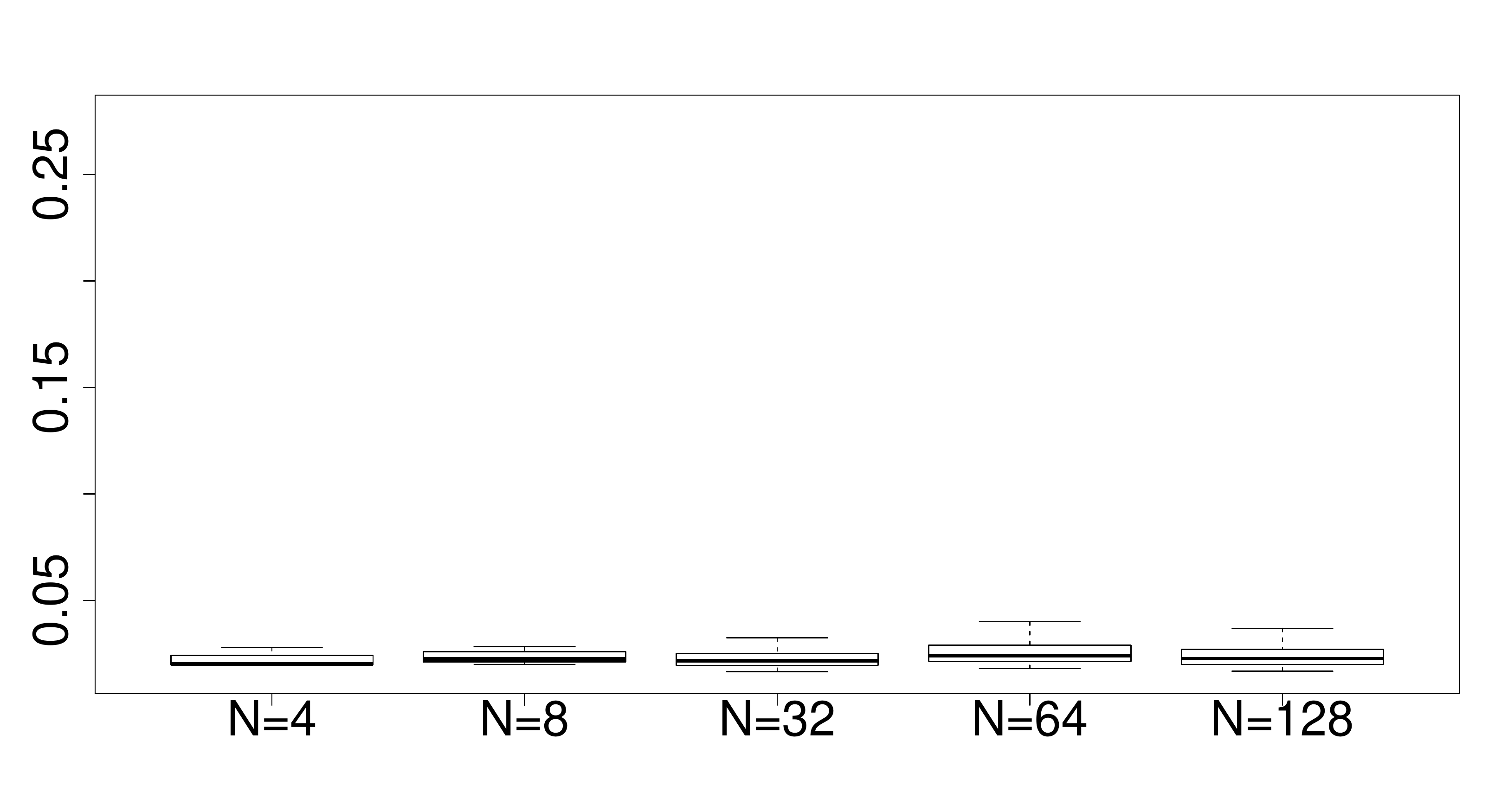}
& \includegraphics[scale=.15]{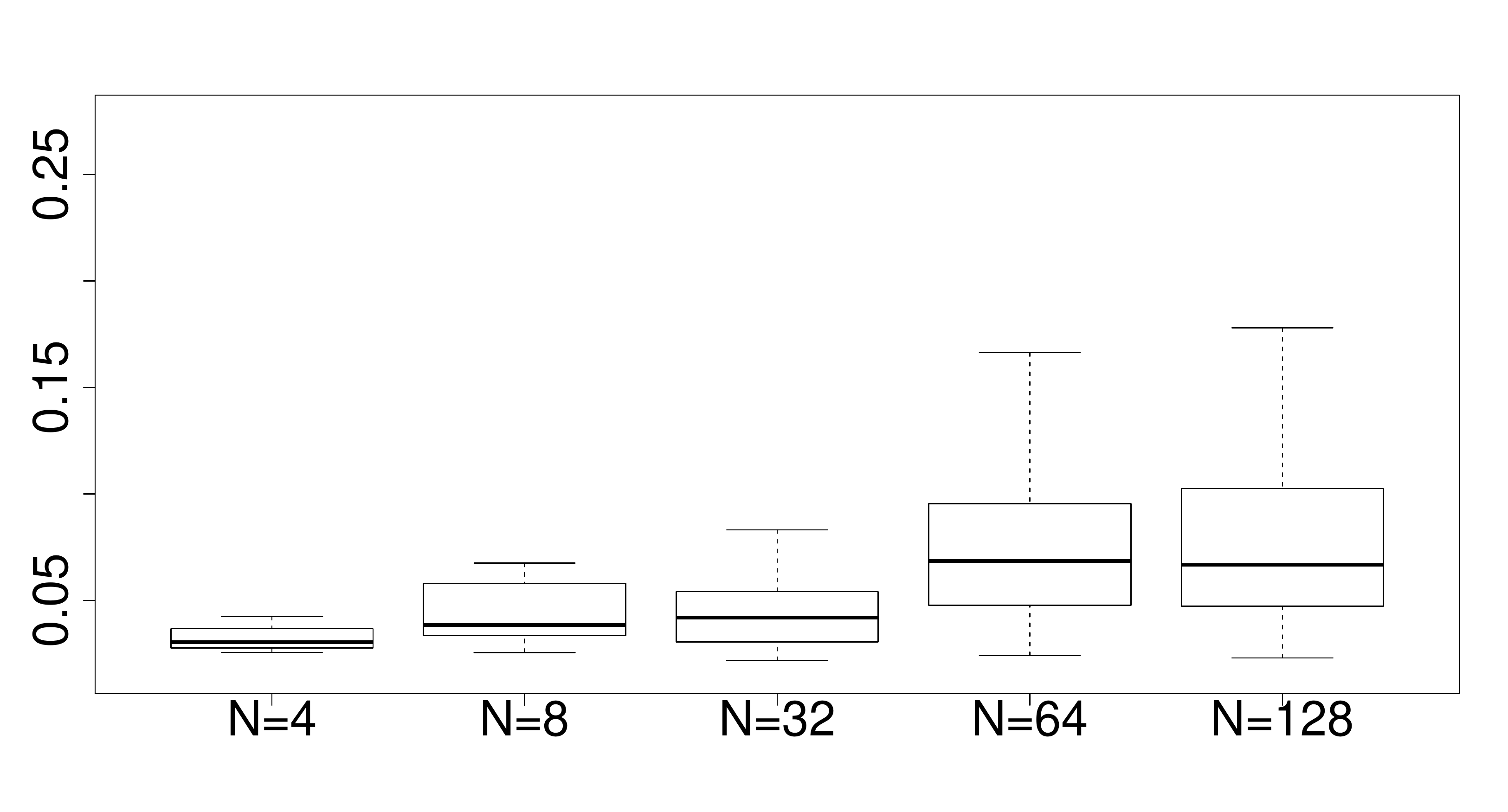}
& \includegraphics[scale=.15]{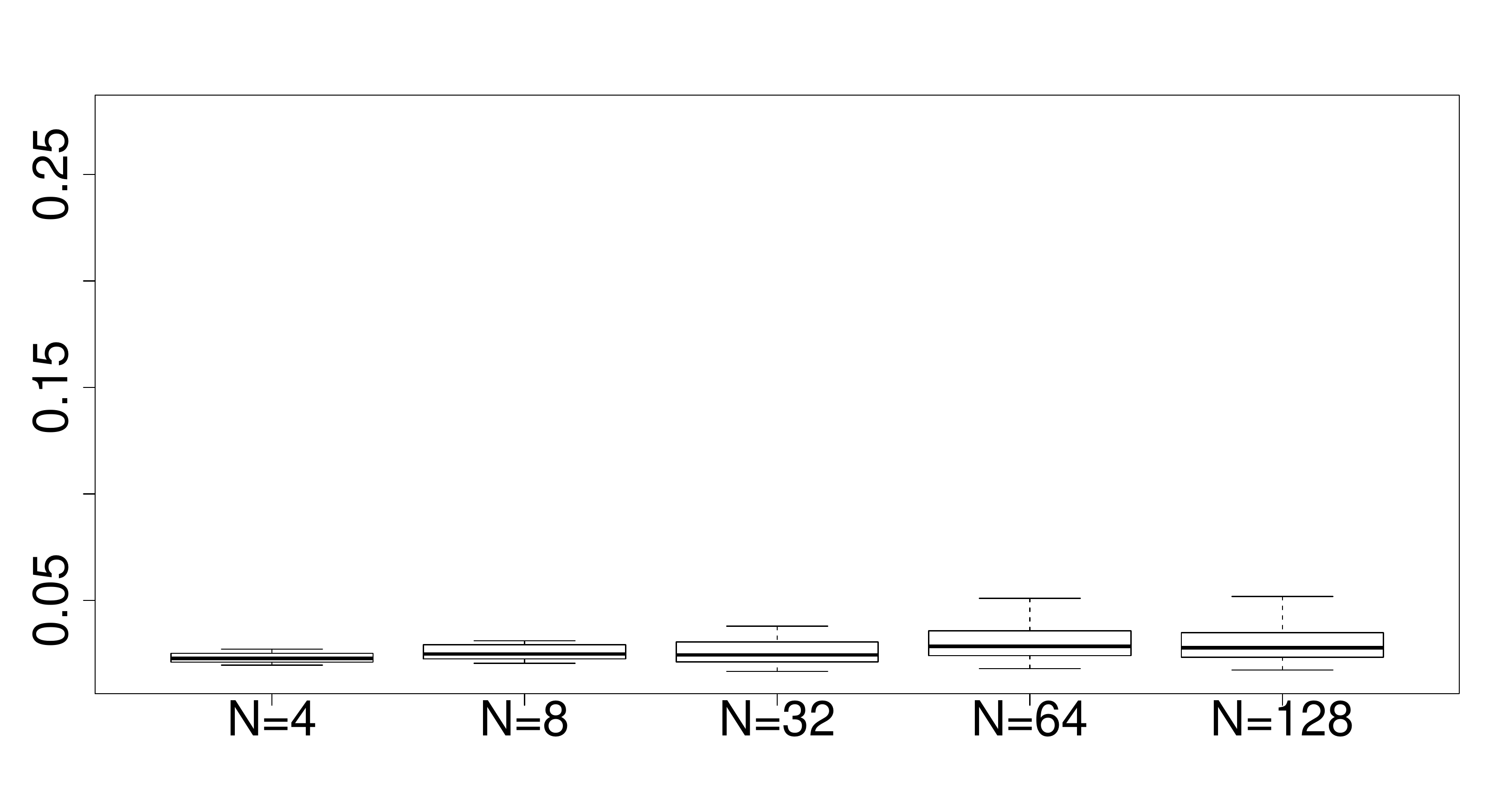} \\
  ES (informed prior) &  VS (informed prior) &  HS (informed prior) \\[2pt]
\hline\\
   \multicolumn{3}{c}{SSIM} \\[2pt]
\includegraphics[scale=.15]{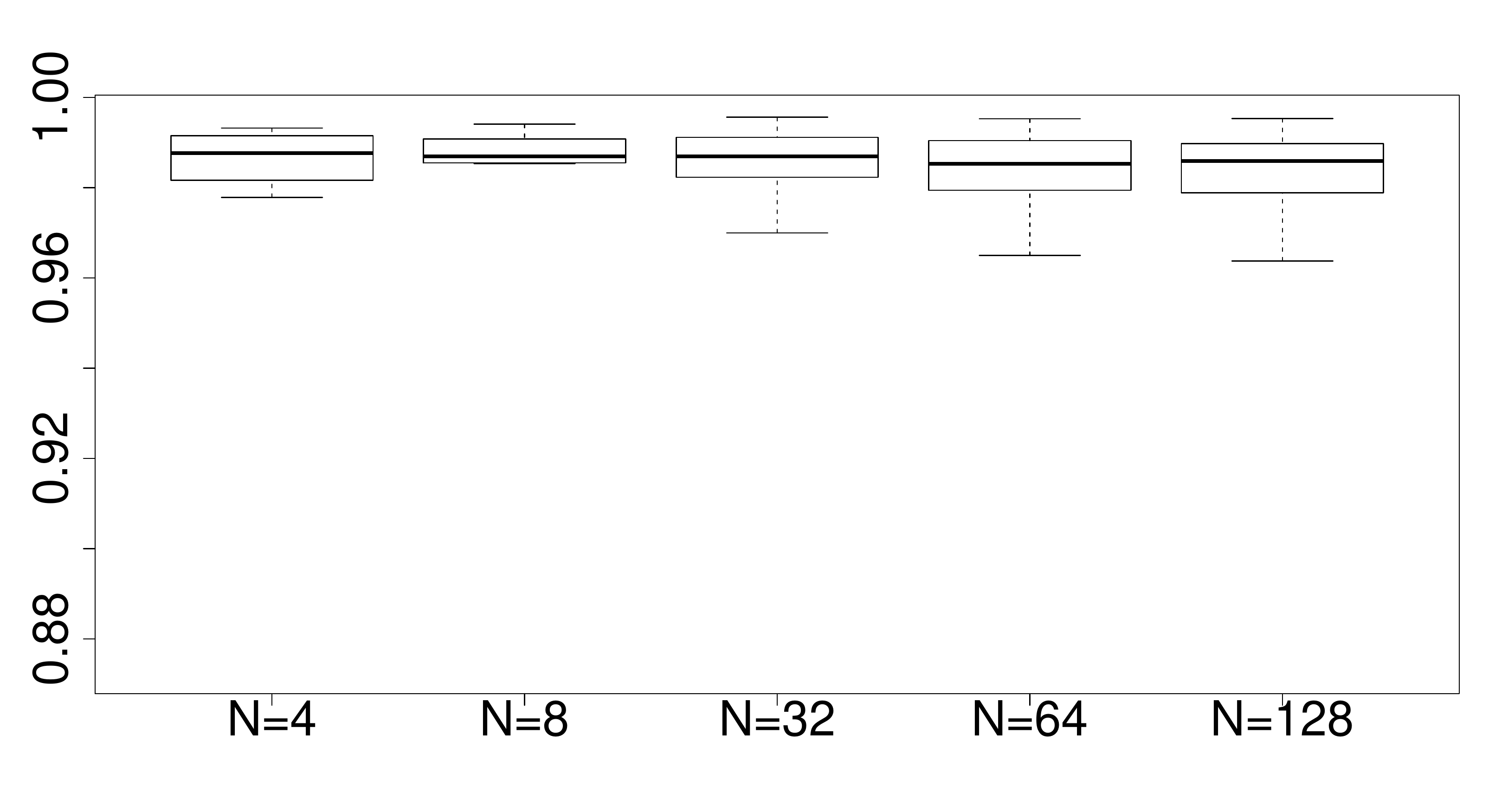}
& \includegraphics[scale=.15]{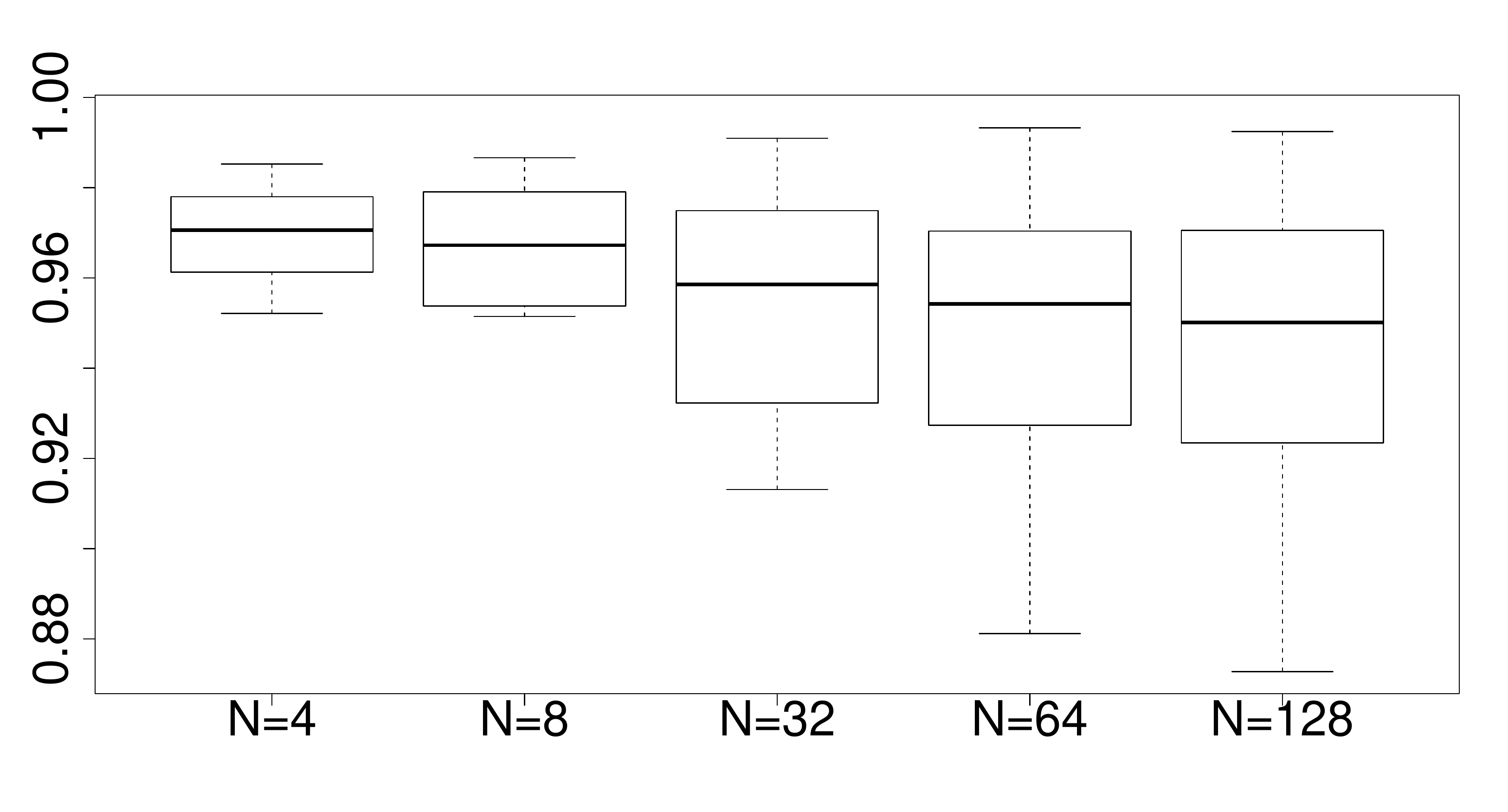}
& \includegraphics[scale=.15]{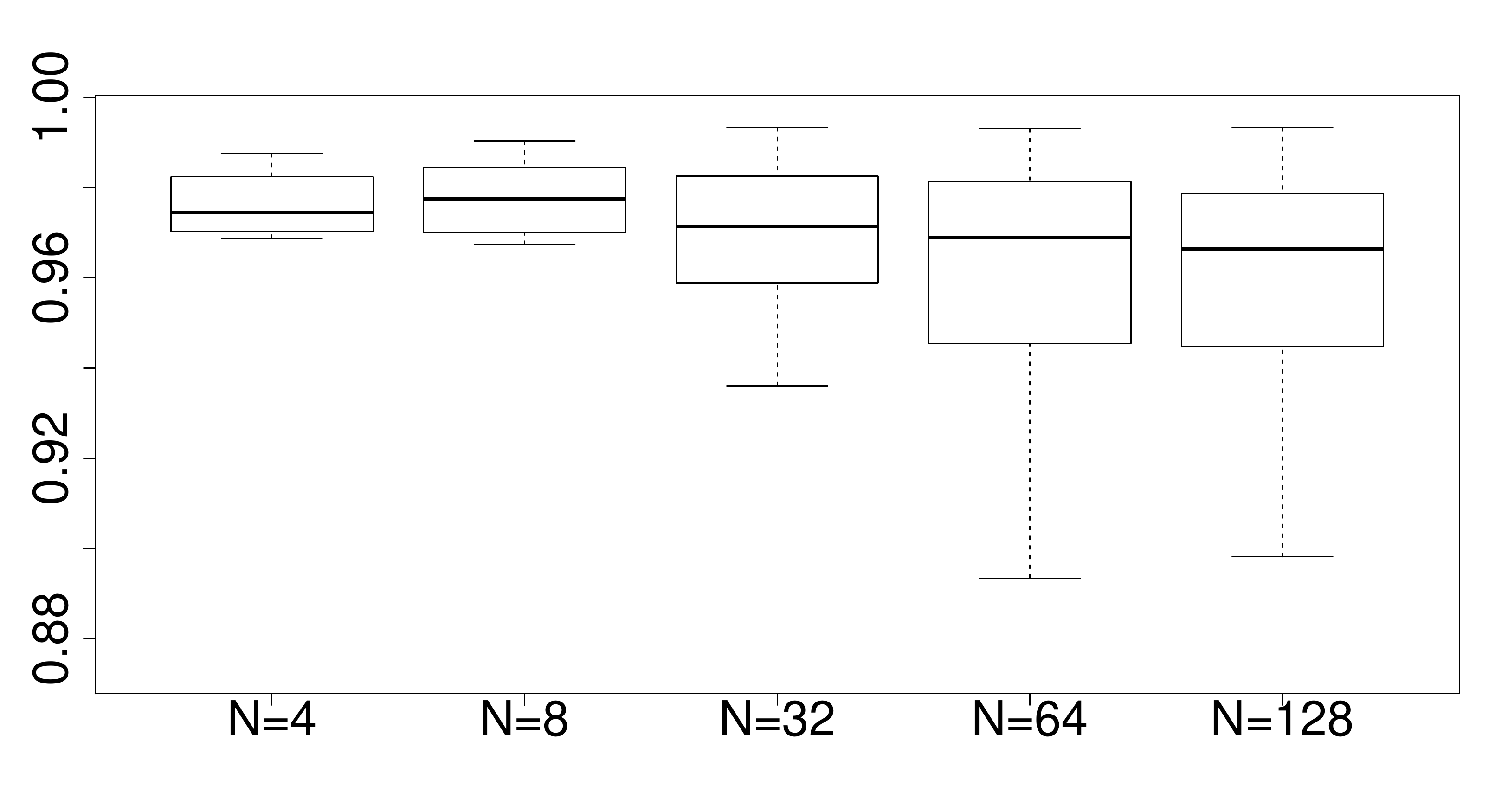}  \\
  ES (standard prior) &  VS (standard prior) &  HS (standard prior) \\
\includegraphics[scale=.15]{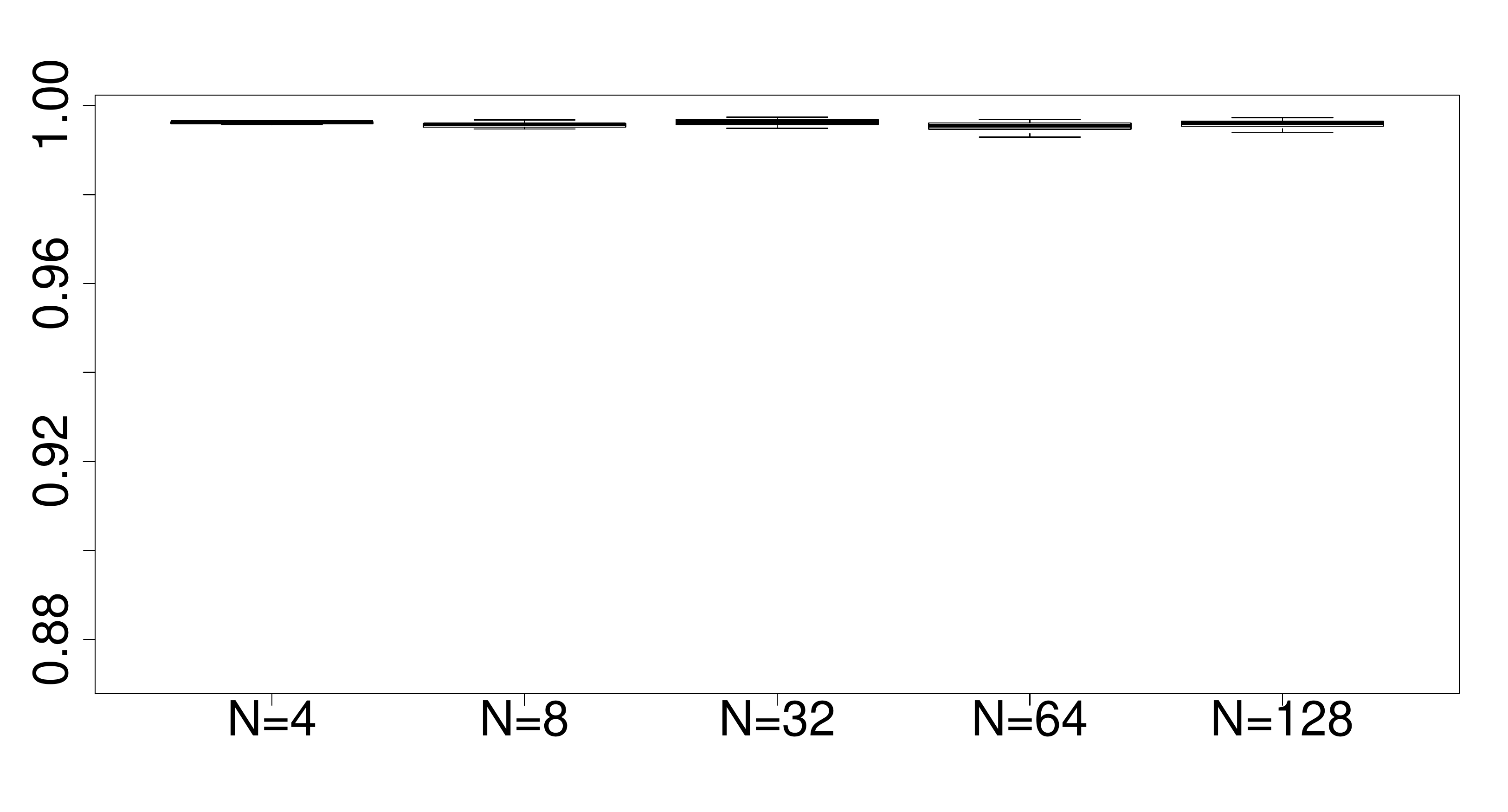}
& \includegraphics[scale=.15]{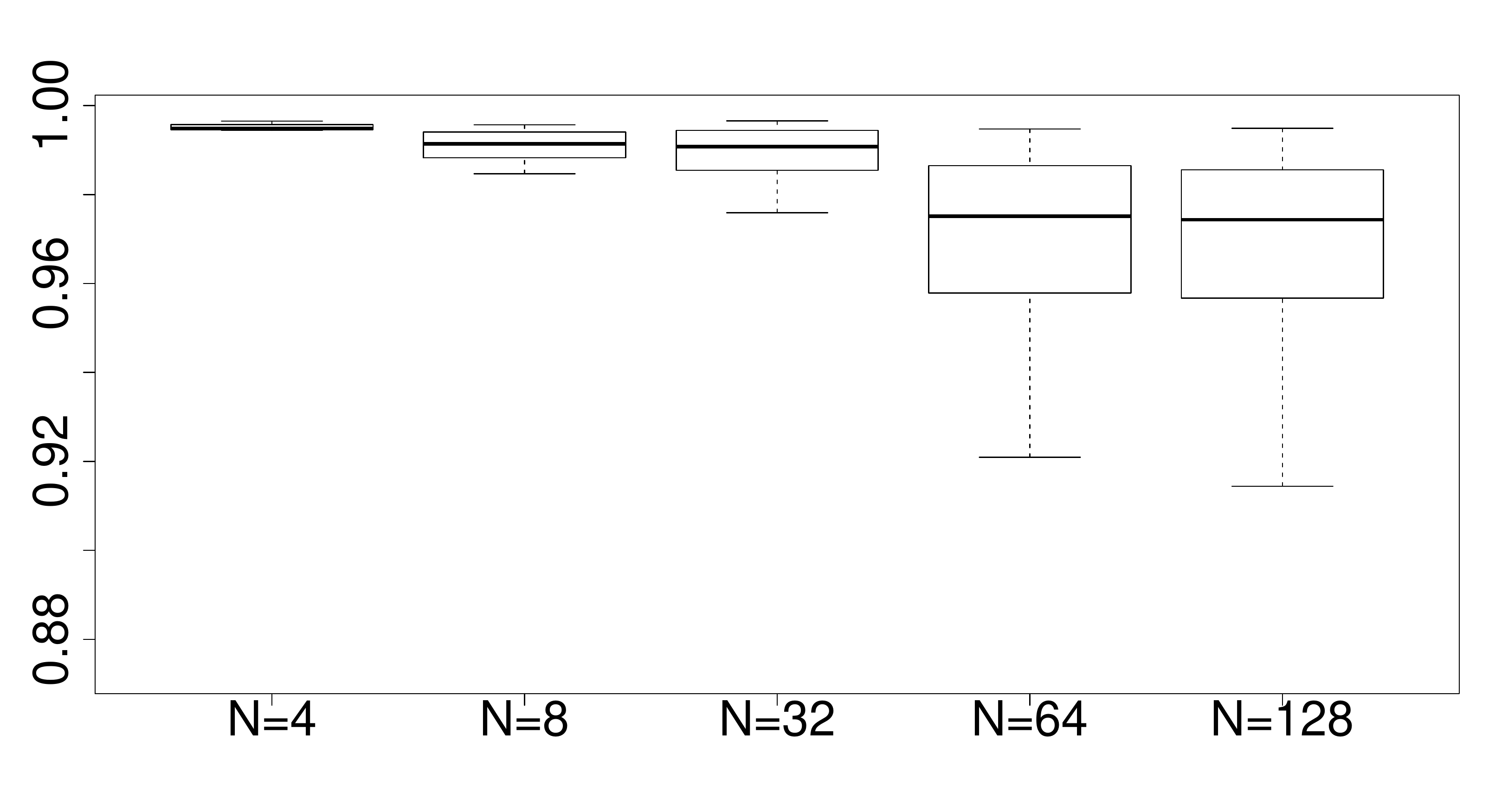}
& \includegraphics[scale=.15]{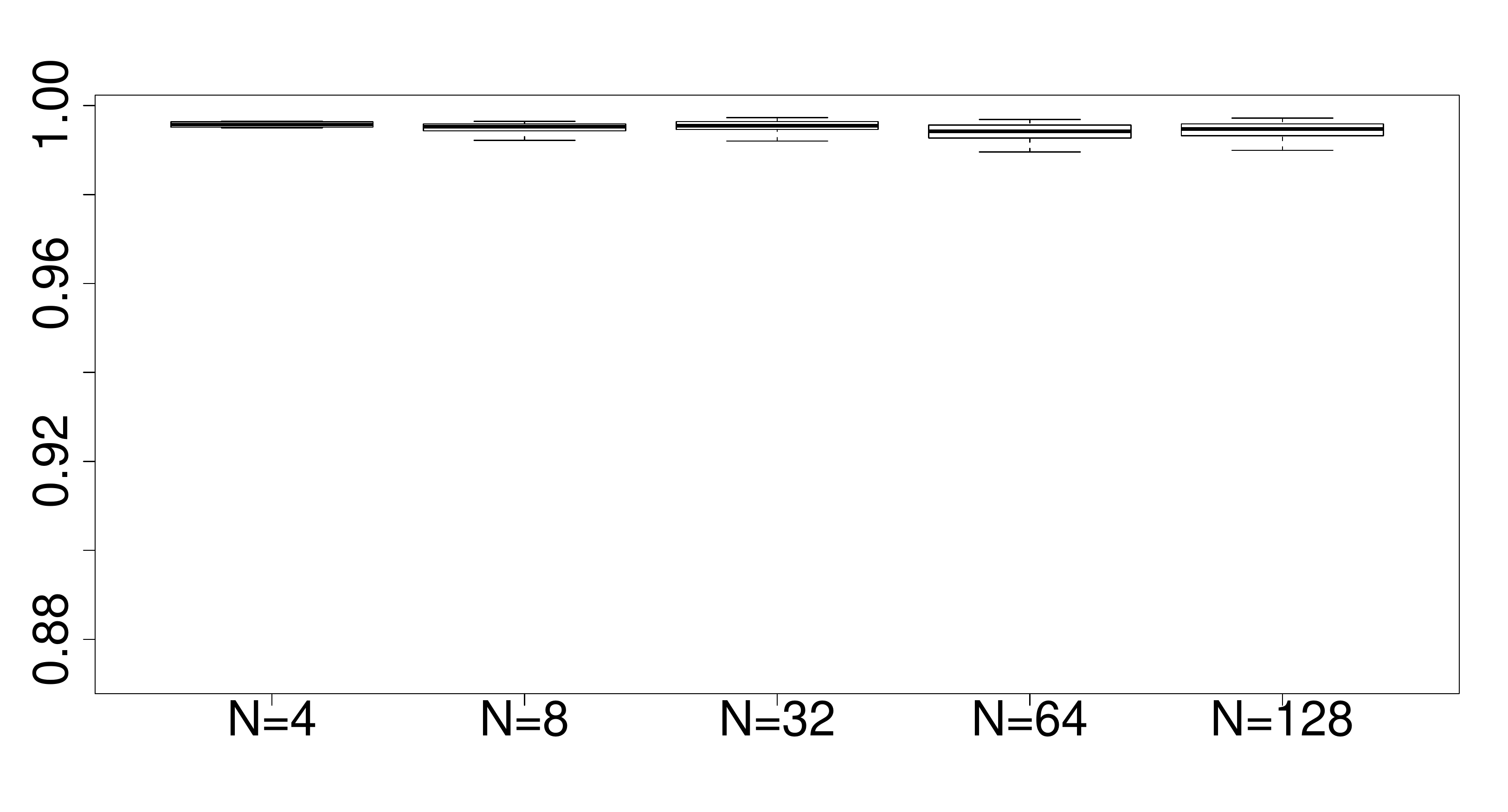} \\
  ES (informed prior) &  VS (informed prior) &  HS (informed prior) \\
\end{tabular}
\caption{RMSE (two upper rows) and SSIM (two lower rows) between $\obsT$ and each sample of $\obs(\ipar)$. Left to right: ES-model; VS-model; HS-model. Models are run with informed  and  standard priors. The number of samples of the right-hand  side varies for each model (N = 4, 8, 32, 64, 128). The RMSE is expected to be as close to 0 as possible, whereas an ideal SSIM would be as close to 1 as possible. }\label{fig:obs:boxplotSSIM}
\end{figure}

In order to investigate the spatial structure of the observables, a metric assessing structural feature, namely, the SSIM, is computed. 
For each generated sample,  in order to assess the overall error between the signals, the SSIM and RMSE are computed between the true and the recovered observable. The values of the metrics are summarized in boxplots in  Figure~\ref{fig:obs:boxplotSSIM}. 
 Comparable to Figure~\ref{fig:obs:PIThisto}, Figure~\ref{fig:obs:boxplotSSIM} shows that the ES-model and HS-model provide
better results than does the VS-model in terms of recovering the
observables $\obs(\ipar)$.  The metrics tend to have more variability for the
VS-model and HS-model when the number of samples increases. This
variability likely comes from the overdispersion of the outputs of
the VS-model.  We note, however, that the overall range of the bulk
of the distribution (box) stays reasonably narrow.  The hybrid score thus appears to be a good compromise between the ES-model and VS-model.

As a conclusion, the overall method shows a wide range of
results in this experimental setup.  In terms of capturing the
parameter field and the intensity and the variability of $m$, the VS- and HS- scores show better
agreement with the true one than does the ES-based model. For
observables, however, the ES-model shows good results in capturing
statistics of the data $\obsT$.  As expected, the informed priors help
capture the parameter $m$ better, as seen on Table
\ref{tab:MAP:metrics} whereas   the standard prior case gives a better calibration between $\obsT$
and $\obs(\ipar)$ and more accurate intensity and structure of $\obs(\ipar)$
(see Figures \ref{fig:obs:PIThisto} and \ref{fig:obs:boxplotSSIM})
arguably by relaxing the parameter constrained through the prior.

\section{Model problem 2: Parameter identification in power grid
  applications governed by DAEs}
\label{sec:num:power:grid}

\begin{figure}[h!]\centering
  \begin{tabular}{cc}
    \includegraphics[width=.30\textwidth]{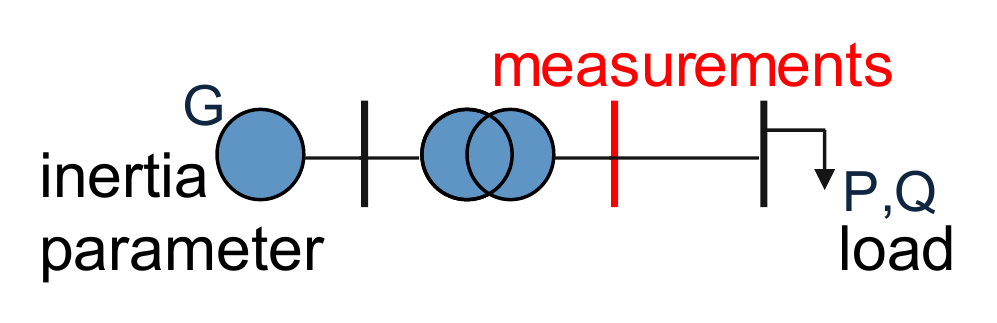}
\end{tabular}
\caption{Power grid diagram.\label{fig:power:grid:diagram}}
\end{figure}

Next we probe the proposed scores on a power grid inverse model
problem governed by an index-1 DAE system. This model 
incorporates an electromagnetic machine, a slack bus, and a stochastic
load as illustrated 
in Figure~\ref{fig:power:grid:diagram}.

We model the power grid using
the generator, current, and network equations~\cite{Sauer_B1998}, namely,
\begin{subequations}
  \label{eq:DAE}
  \begin{align}
    \begin{cases}
      \dot{x} &= f(x,y;\ipar) \\
      0 &= g(x,y;\xi) 
    \end{cases}
    \Rightarrow
    F(u,\ipar;\NoiseForcing)=0 ~ \textnormal{a.s.}\,, ~
    u=[x,y]^T\,.
  \end{align}
\end{subequations}
Here $x$ is associated mainly with generators, and $y$ represents
current (part of the generator equations) and the network equations
(Kirchhoff). The full set of equations is given in Appendix
\ref{sec:power:grid:eq}. The unknown (or inversion) parameter here is
$\ipar$ and represents the generator inertia, which is one
dimensional in this example. This parameter can be
interpreted as how fast the generator reacts to fluctuations in the
network. In our case $\xi=[P,Q]^T$
represents fluctuations in the load, where $P$ and $Q$ represent the
real and imaginary resistive components, respectively.

In recent work we have explored estimating the inertia parameters in a
standard 9-bus system given a single known disturbance in the load
from synthetic bus voltage observations~\cite{Petra_2016a}. In this
study we pose the problem as having a small signal disturbance in the
load, which is a discrete process in time. This problem now describes
a realistic behavior of small-scale consumers drawing power from the
grid in an unobserved fashion. 
What we consider to be known is the
distribution of the probabilistic load process. Moreover, we assume
that we measure the power flow (or voltage) at one of the buses (see
Figure~\ref{fig:power:grid:diagram}).

\begin{figure}[h!]\centering
    \subfigure[{$V_r$ simulated}]{\includegraphics[width=.35\textwidth]{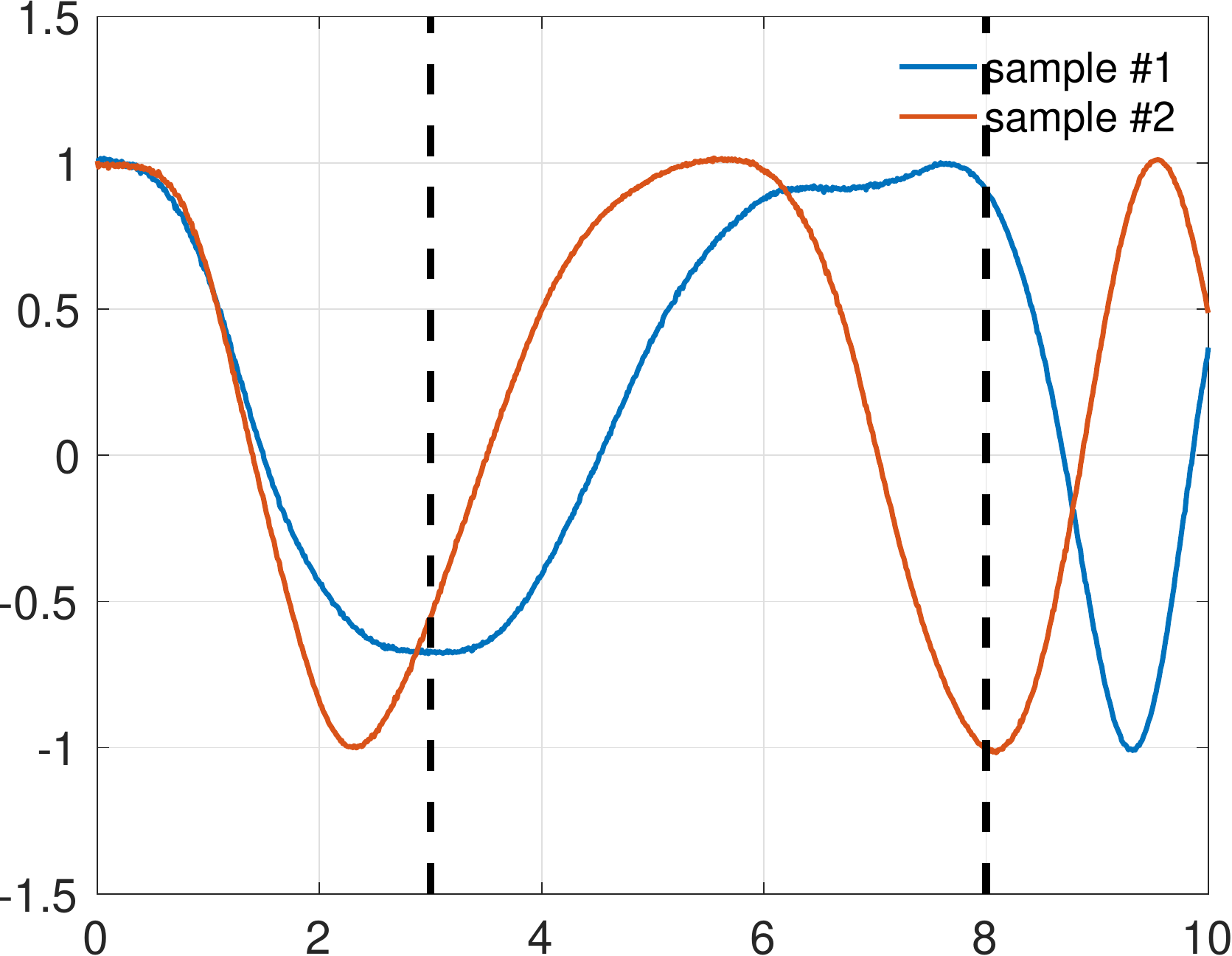}\label{fig:Vr}}
    \subfigure[{$V_i$ simulated}]{\includegraphics[width=.35\textwidth]{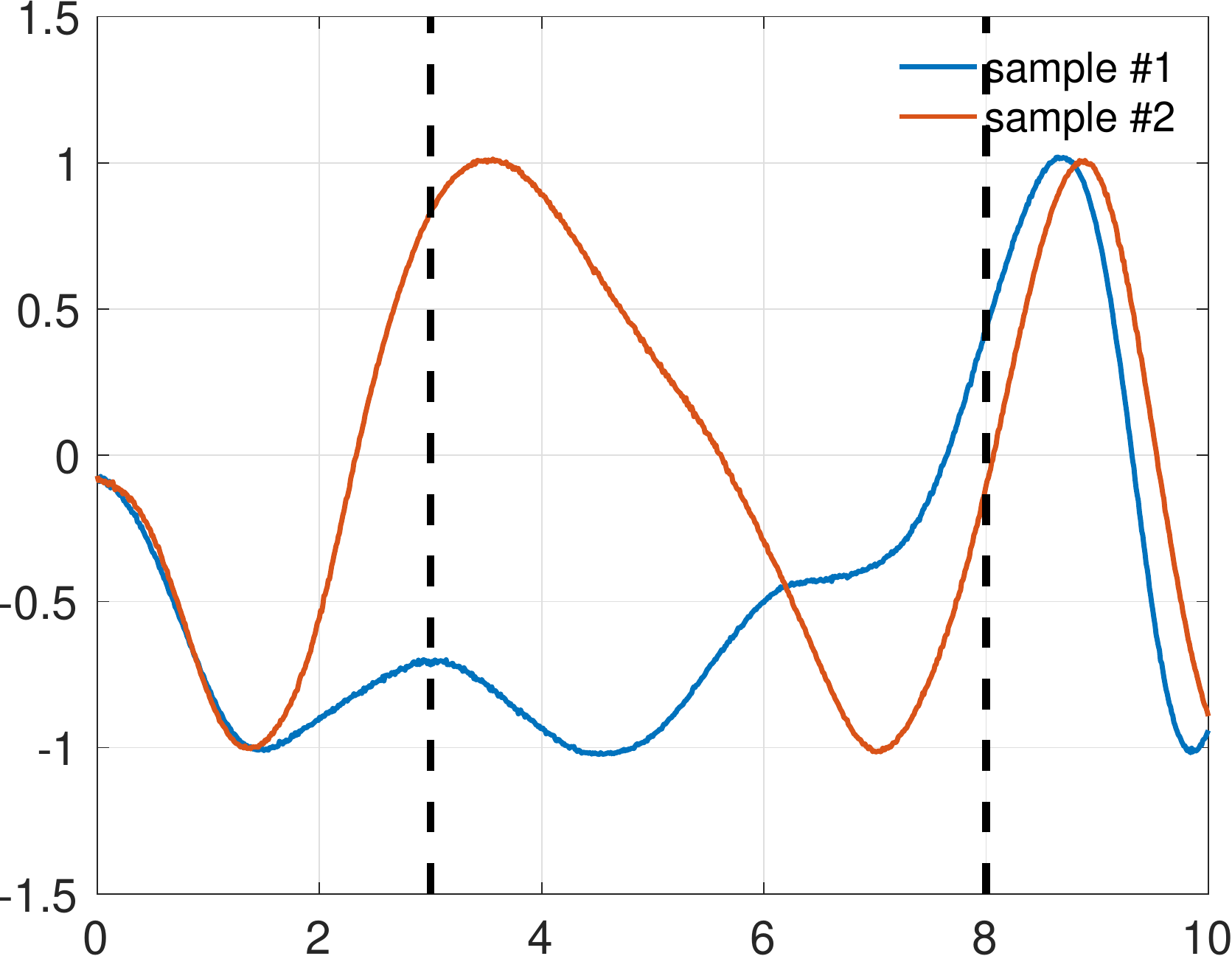}\label{fig:Vi}}
    \subfigure[{$P$ time series}]{\includegraphics[width=.35\textwidth]{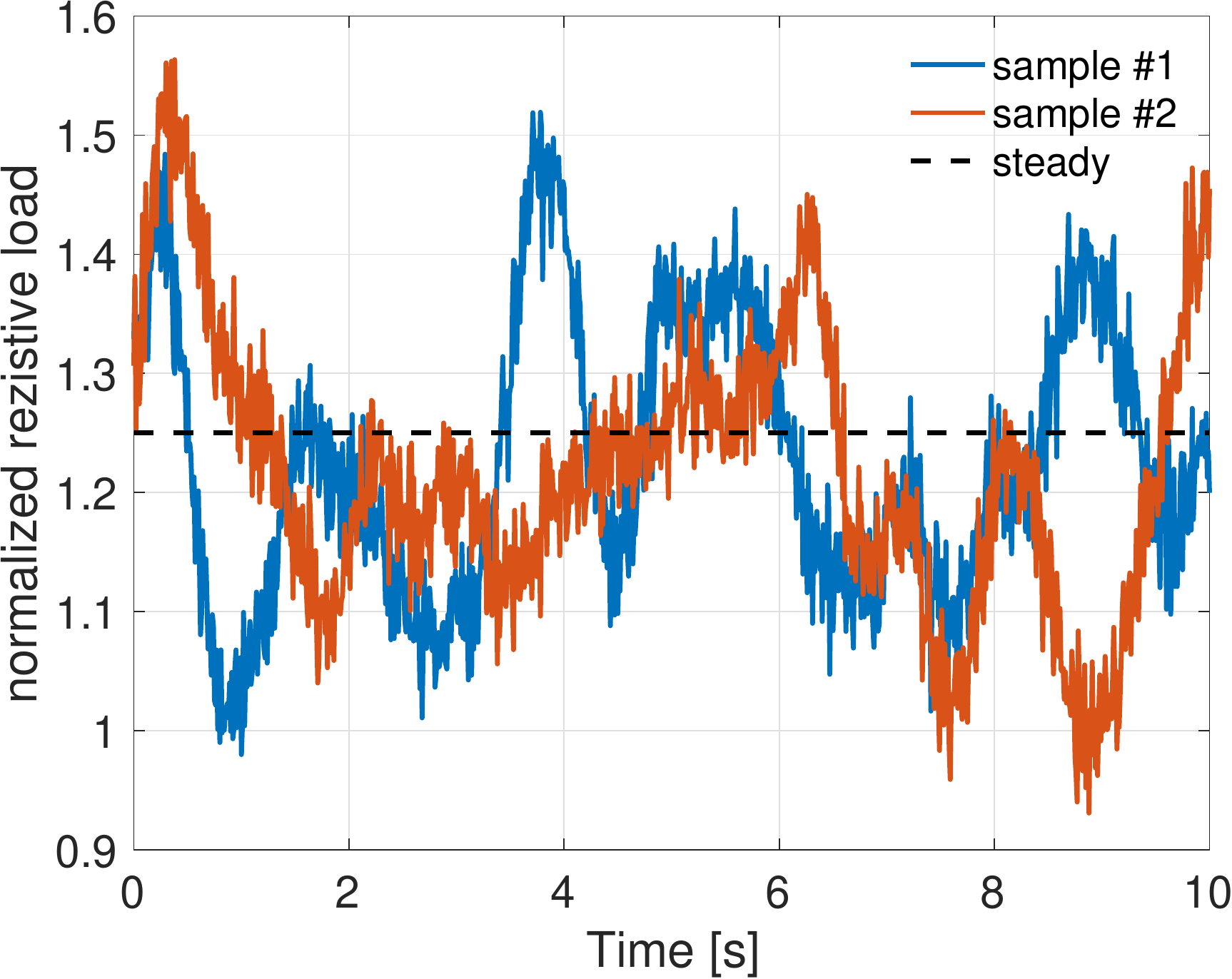}\label{fig:P}}
    \subfigure[{$Q$ time series}]{\includegraphics[width=.35\textwidth]{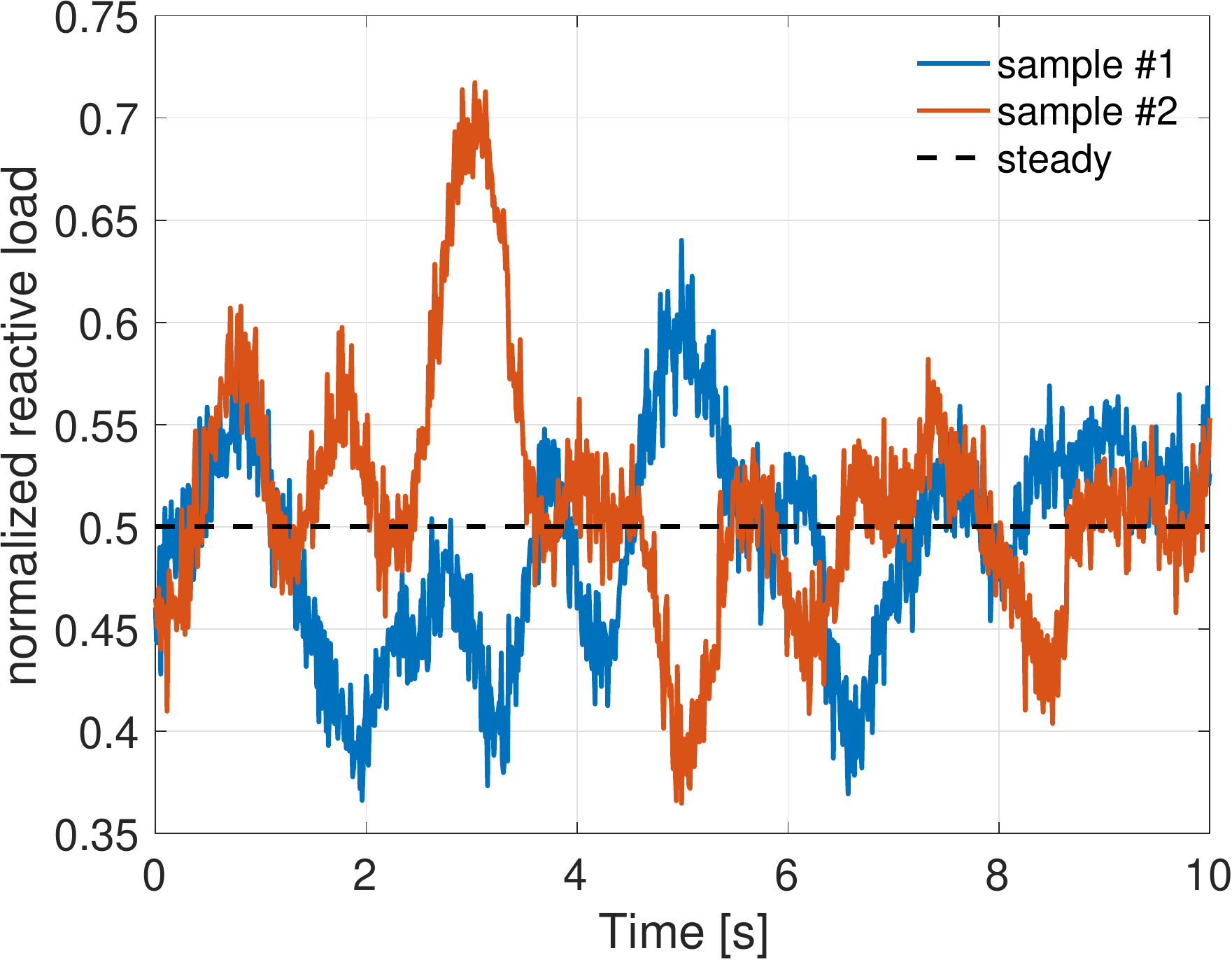}\label{fig:Q}}
\caption{Voltage with (a) real and (b) imaginary parts at the measurement
  location for two samples. Load noise around the stationarity
  baseline for the (c) resitive and (d) reactive components. \label{fig:VrVi:PQ}}
\end{figure}

\subsection{Computational experiment}

We start with the system at dynamic steady state having
$\frac{\partial x}{\partial t} = \frac{\partial y}{\partial t} = 0$
and $\xi = \left[\overline{P},\overline{Q}\right]^T =
[1.25,0.5]^T$. The system is integrated with backward Euler with
stochastic forcing. This discretization is then equivalent to
first-order weak convergence in the mean square error sense of the
stochastic DAE. Higher-order methods have been used as well, without
noticing any qualitative differences for our current setup. The time series generated by the
stochastic forcing is given by two independent stationary processes
with known distribution $\pinoise$:
\begin{align}
  \label{eq:pg:Noise}
  \xi=[P,Q]^T, P(t) \sim {\cal N}\left( \overline{P}, 0.1^2
  \,\CovFcn(h) \, 
  \right),~
Q(t) \sim  {\cal N}\left( \overline{Q}, 0.05^2 \,  \CovFcn(h) \,
\right),~ \CovFcn(h)=e^{-\frac{h^2}{0.002}} +
0.1,
  \end{align}
where $h=t-t'$ and, as before, 
$\mathrm{Cov}(P(t),P(t'))=\CovFcn(h)$. Some realizations of these time
series are shown in Figure~\ref{fig:P}-\ref{fig:Q}.

The simulation time is $T=10$ seconds, with a time step of $\dt = 10^{-2}$ (10
ms). From the 10-second window we extract 5 seconds (seconds 3 to 8) to avoid
initialization or mixture issues. We consider an ensemble of 1,000
samples integrated with this time step, with $\NEns=800$ being considered as
numerical simulations and $n=1,2,\dots,200$ set aside for observations. In this
experiment we do not consider observational noise ($\varepsilon_{\rm obs} \equiv 0$) and do not need to
use any regularization; this is equivalent to an uninformative or flat
prior. The optimization problem becomes a univariate
unconstrained program, and therefore we approximate the gradient with
finite differences. Nevertheless, one can compute the gradients via
adjoints as has been done for model problem 1 in Section~\ref{sec:poisson} as
well.

\subsection{Analysis of the results}

The time-dependent setting allows us to analyze different aspects of
the inverse problem solution.
For instance, in the steady-state
case such as the first model problem, the spatial domain is
fixed, and in general so is the number of observations. In the unsteady
example, one typically has control over the observation window. To
this end, we begin by exploring the effect of adding observations to
the inference process and thus using the mean score
\eqref{eq:mean:score}. Specifically, we compute the score values at
integer values of the parameter, from 1 to 35, and use 1 to 200
batches of observations or validation samples. In other words we explore the mean score
values $\ScoreFcn_n$ with $n=1,\dots,200$; i.e., $\obsT^{(1,\dots,n)}=[\obsT^{(1)},\obsT^{(2)},\dots
  \obsT^{(n)}]^\top \in  \mathbb{R}^{n \times \Nobs}$. One batch is the time
series obtained under a stochastic forcing realization. Each batch is
the result of a 5-second simulation, and this assumes that the
distribution is stationary for the entire inference window; i.e., the
distributions do not change over time. The
observations of the voltage $x_{11}$ and $x_{14}$ corresponding to the
middle bus are taken at every time step.

\begin{figure}[h!]\centering
  \begin{tabular}{cc}
    \includegraphics[width=.45\textwidth]{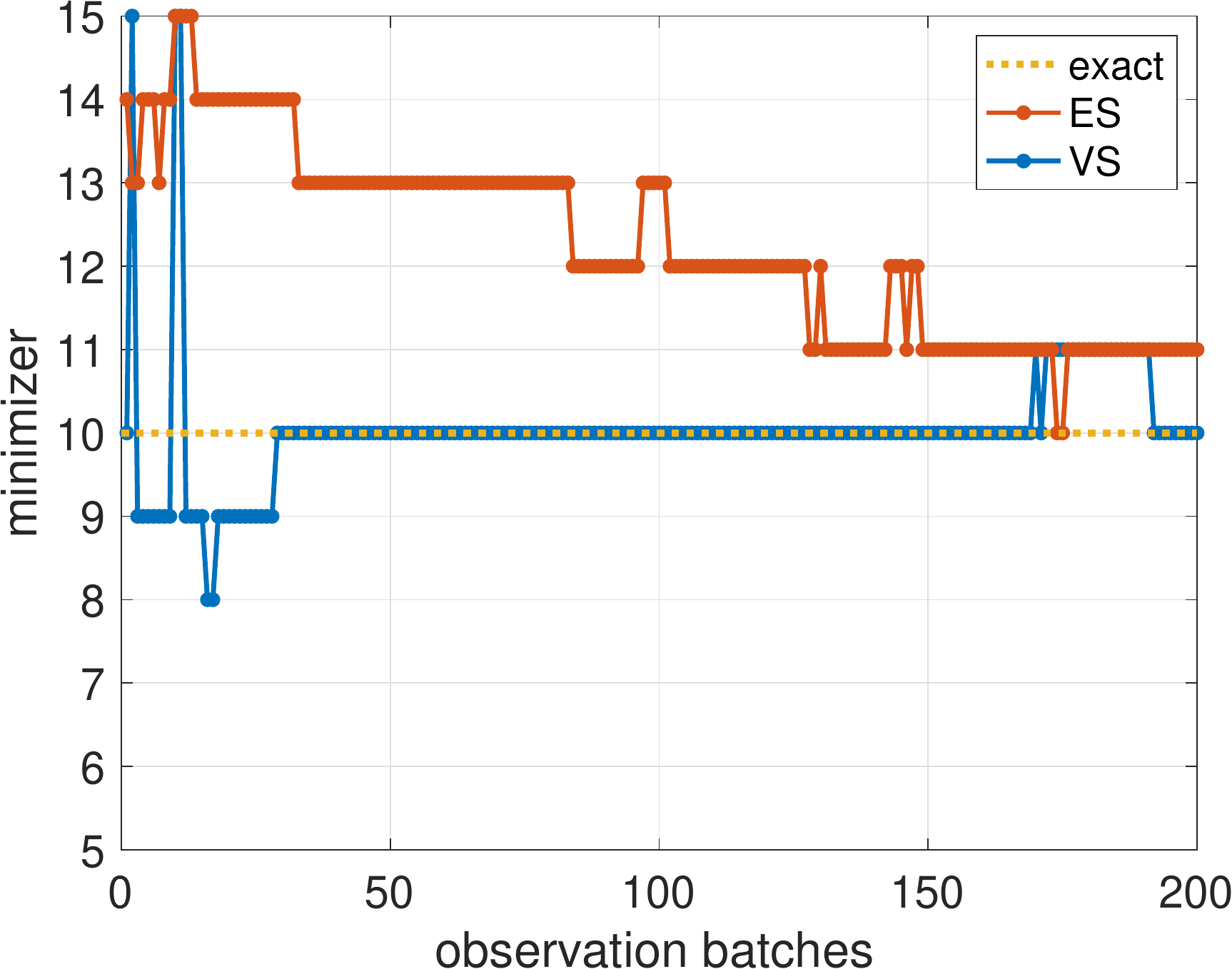}&
    \includegraphics[width=.45\textwidth]{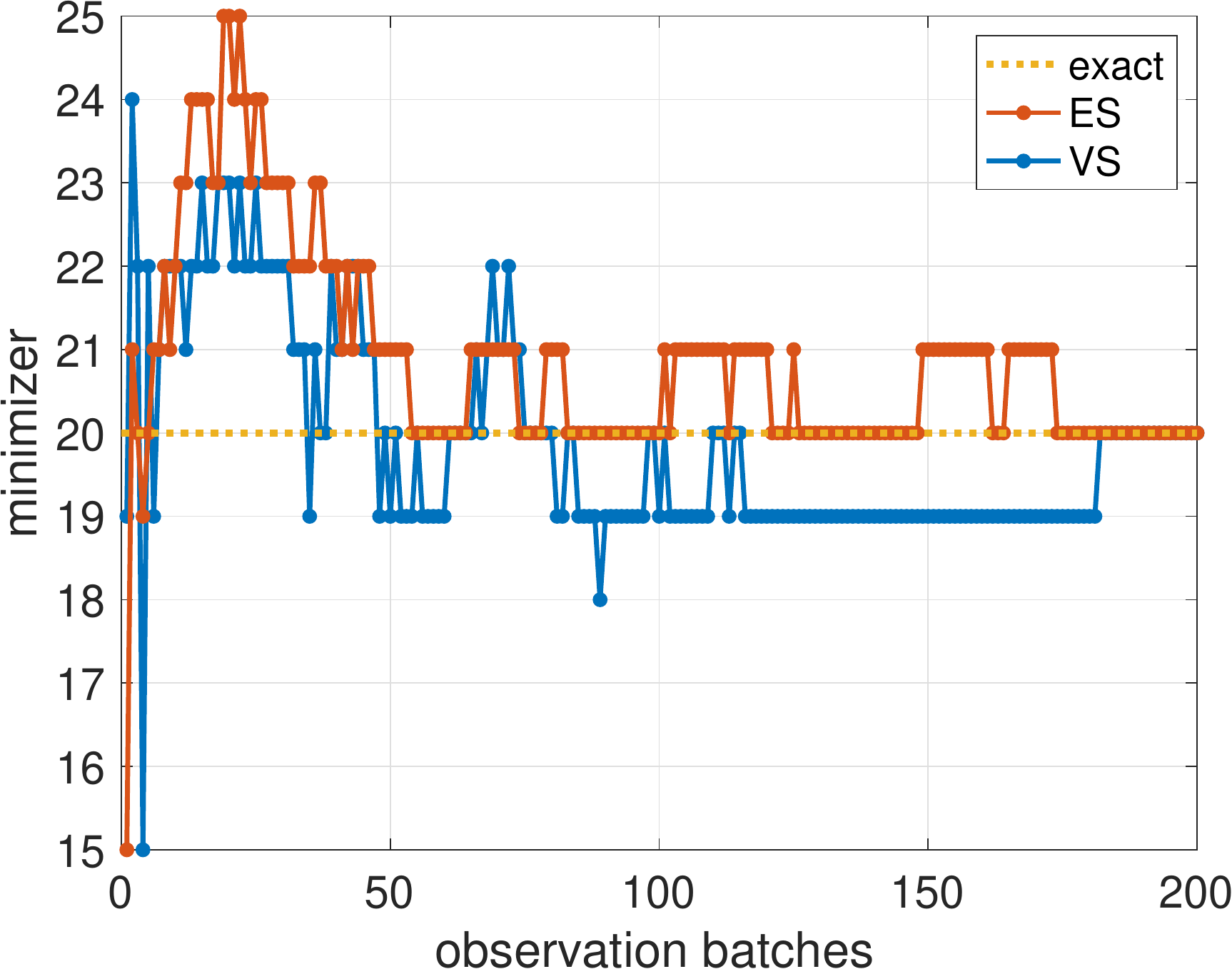}
\end{tabular}
\caption{Grid search of the minimizer as a function of observation
  batches. The exact parameter value is 10 for the left panel and 20 for
  right one. As the number of observations increases, the minimizer
  converges to the true value. \label{fig:optim:convergence}}
\end{figure}

We illustrate the results in Figure~\ref{fig:optim:convergence} for
the ES-model and VS-model computed with respect to two exact
values of the parameter, 10 and 20, respectively. We observe that both the energy
and variogram scores converge to the exact value, with the variogram
score converging much faster than the energy score especially when the
exact value of the parameter is 10. We also remark that convergence
guarantees are not easy to ascertain a priori; as reflected in the
figure, a different number of observations are necessary in order to
reach an accurate conclusion. One possible strategy to mitigate this
issue is to use two different scores and observe the system until they
are in agreement and  do not change with additional observations.

\begin{figure}[h!]\centering
  \begin{tabular}{cc}
    \includegraphics[width=.45\textwidth]{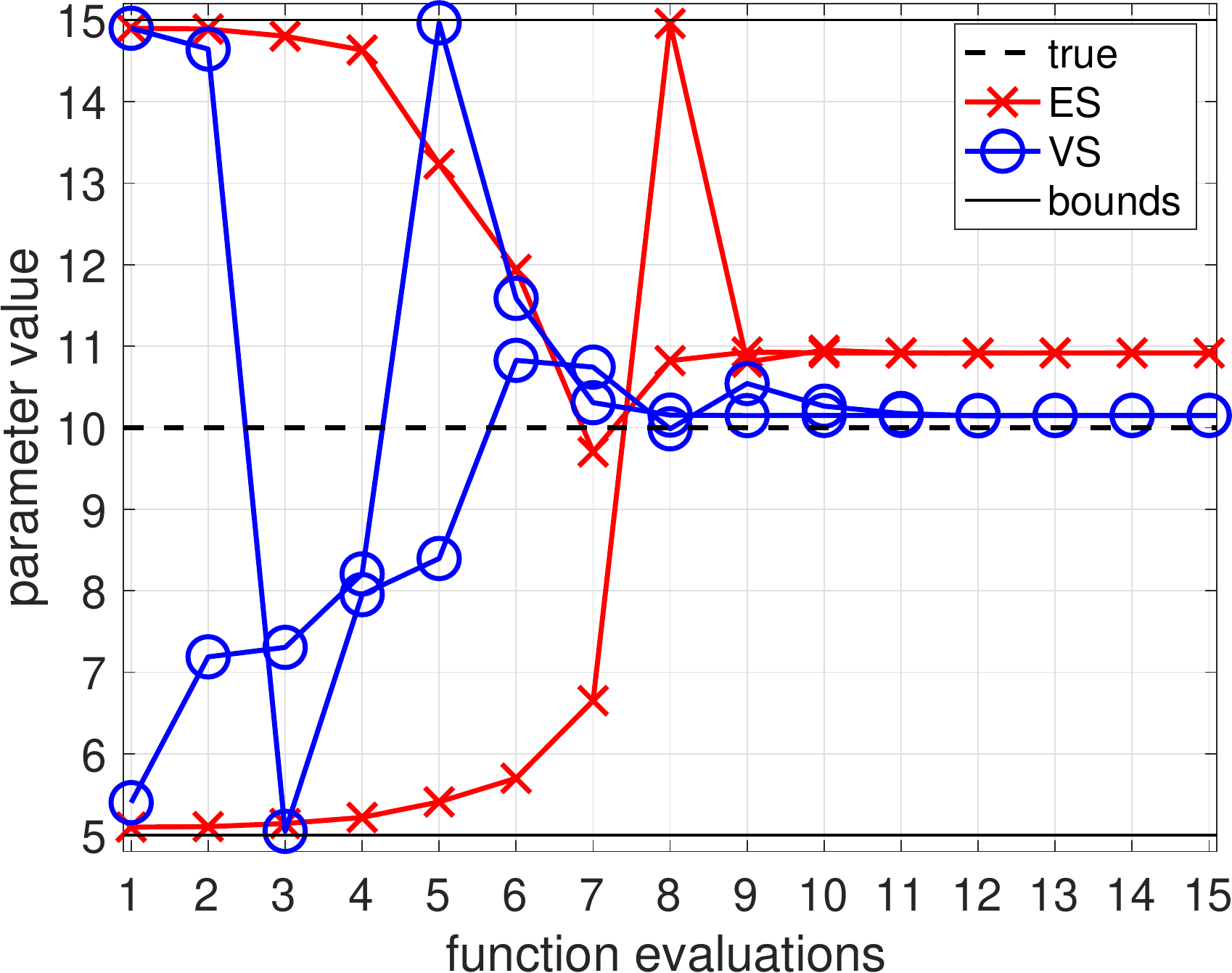}
    &
    \includegraphics[width=.45\textwidth]{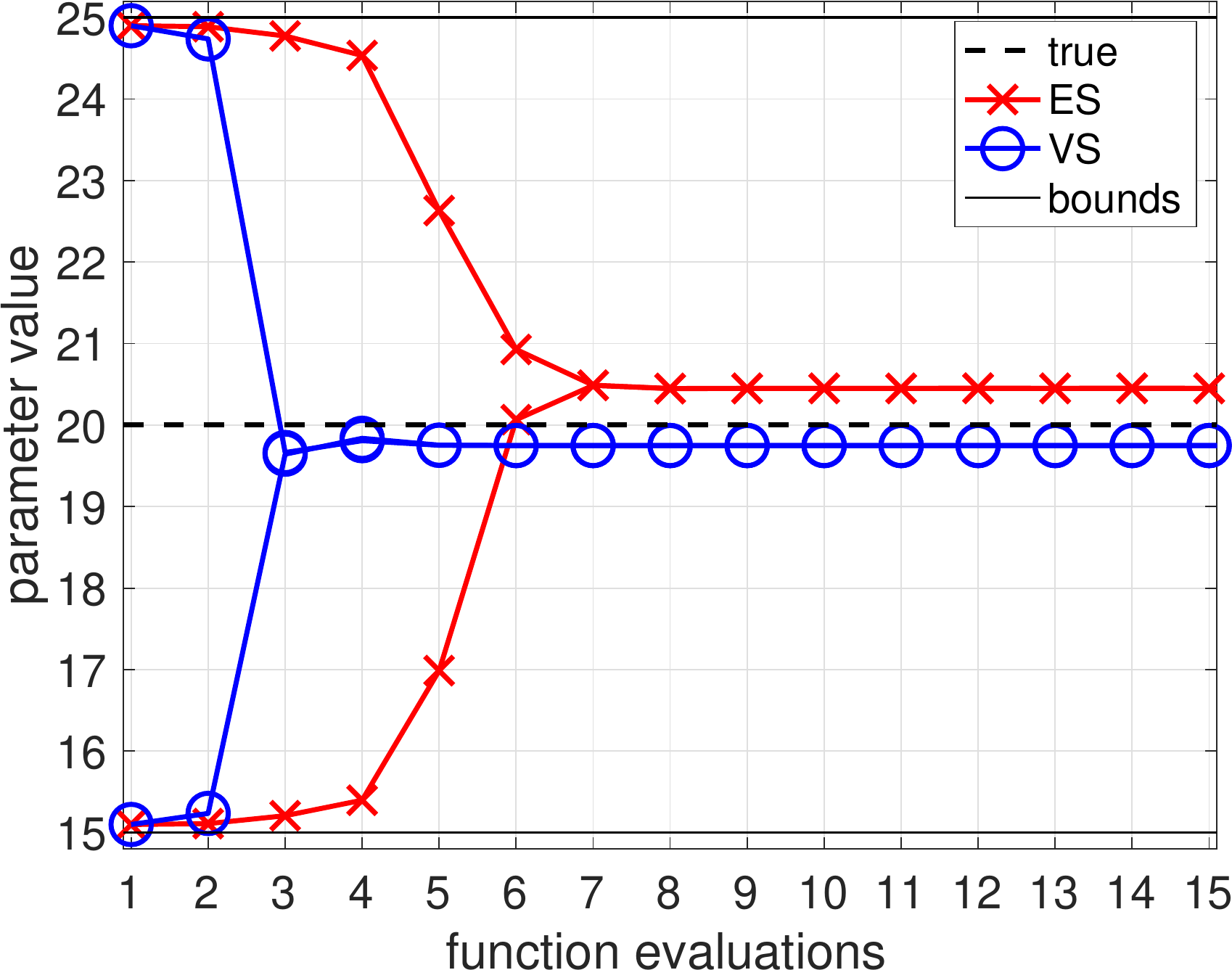}
\end{tabular}
\caption{Reconstructed parameter for $\ipart = 10$ (left) and
  $\ipart = 20$ (right). The results obtained with the energy score
  are shown in red (cross) and for the variogram score in blue
  (circles). The bounds
  $\ipart\,\pm \,5$ used in the optimization solver are shown
  with a black solid 
  line. Two experiments are carried out for each score with the
  initial guess being the high and low bound values. \label{fig:optim:minimization}}
\end{figure}
The numerical results are carried out in Matlab using the default 
optimization solver. In Figure~\ref{fig:optim:minimization} we show
the reconstructed parameter values as a function of function
evaluations for the exact parameter values (black dashed line) 10
(left) and 20 (right). The results for the energy score are shown in
red (cross) and for the variogram score in blue (circles). The bounds
set to truth$\,\pm\, 5$ used in the optimization solver are shown with a black
solid line. This figure shows that with both scores the optimization
solver converges to a relatively good estimate of the exact parameter
value in a relatively small number of function evaluations.

\section{Conclusions\label{sec:conclusion}}

We have presented a statistical treatment of inverse problems governed
by physics-based models with stochastic inputs. The goal of this study
is to quantify the quality of the inverted parameter field, measured
by a comparison with the ``truth,'' and the quality of the recovered
observable: for example,  given the inverted parameter field, quantify how
well we fit the distribution of the observable. The end goal of our
study is to introduce an inverse problem formulation that
facilitates the integration of data with physics-based models in order
to quantify the uncertainties in model predictions. To this end, inspired
from statistics, we propose to replace the traditional least-squares
minimization problem---which minimizes the norm of the misfit between
data and observables---with a set of loss functions that describes
quantitatively the distance between the distributions of model
generated data and the distribution of observational data. We refer to 
these metrics as scores, as known in the statistics community.

To compute the maximum utility or a posteriori point for the proposed
inverse problem, we solve an optimization problem constrained
by the physics-based models under stochastic inputs with a
quasi-Newton limited-memory algorithm. For efficient calculation of
the gradient of the objective with respect to the inversion
parameters, we derive adjoint-based expressions.  Several
challenges are associated with solving such optimization problems. First,
these inverse problems are large scale, stemming from discretization
of the parameter field in the case of PDE-based models or the size of
the power grid network. Second, although we employ an efficient method
to calculate derivatives, the number of adjoint solves increases with the number of samples and are coupled. As we indicated above, however,
the communication during the adjoint calculations follows a fixed
pattern and can be optimized for, arguably resulting in overall scalable strategies. Third, structured error in measurements
(data) requires special attention, 
and convergence guarantees are not easy to
ascertain a priori. Nevertheless, as we illustrate in the second
model (see Fig. \ref{fig:optim:convergence}), one can use multiple sets of data to ascertain and
mitigate potential convergence issues. 

We have studied the performance of the proposed formulation in the
context of two applications: a coefficient field inversion for
subsurface flow governed by an elliptic PDE with a stochastic source
and a parameter inversion for power grid governed by DAEs. In both
cases the goal was to obtain predictive probabilistic models that
explain the data.

\section*{Acknowledgments}
We thank Michael Scheuerer for providing  helpful comments on
scoring functions and Umberto Villa for helpful discussions about hIPPYlib. The work of C. G. Petra was performed under the auspices of the U.S. Department of Energy by Lawrence Livermore National
Laboratory under Contract DE-AC52-07NA27344. This material also was based upon work supported by the U.S. Department of Energy, Office of Science, under contract DE-AC02-06CH11357. 

\bibliographystyle{siamplain}

\vfill
\begin{flushright}
  \scriptsize \framebox{\parbox{3.2in}{
The submitted manuscript has been created by UChicago Argonne, LLC,
Operator of Argonne National Laboratory (“Argonne”). Argonne, a
U.S. Department of Energy Office of Science laboratory, is operated
under Contract No. DE-AC02-06CH11357. The U.S. Government retains for
itself, and others acting on its behalf, a paid-up nonexclusive,
irrevocable worldwide license in said article to reproduce, prepare
derivative works, distribute copies to the public, and perform
publicly and display publicly, by or on behalf of the Government.  The
Department of Energy will provide public access to these results of
federally sponsored research in accordance with the DOE Public Access
Plan. \url{http://energy.gov/downloads/doe-public-access-plan}.
}}
\normalsize 
\end{flushright}
\clearpage
\newpage
\appendix

\section{Power grid equations\label{sec:power:grid:eq}}

Below we present the equations extracted from \cite{Sauer_B1998} for the power grid example discussed
in Sec. \ref{sec:num:power:grid}. The one generator 3-bus system is
described by index-1 DAEs. Here we have 7 
differential equations and 8 algebraic equations. The differential
variables are the first seven variables, with the rest being algebraic.
\begin{footnotesize}
  \begin{align*}
\dot{x}_1=&-376.99111843077515 + x_2\\
\frac{\ipar}{23.64} \dot{x}_2=&47.70113037725341 - 0.09968102073365231 x_2 -
7.974481658692184 (x_4 x_{8} + x_3 x_{9} + 0.0361 x_{8} x_{9})\\
\dot{x}_3=&0.11160714285714285(x_5- x_3) -
0.009508928571428571 x_{8}\\
\dot{x}_4=&-3.2258064516129035 x_4 + 1.0938709677419356 x_{9}\\
\dot{x}_5=&-0.012420382165605096 \exp(1.555 x_5) + 3.1847133757961785 (x_{7} - x_5)\\
\dot{x}_6=&0.5142857142857145 x_5 - 2.857142857142857 x_{6}\\
\dot{x}_7=&109.644151839917 - 18 x_5 + 100 x_{6} - 5 x_{7} - 100 \sqrt{x_{10}^2 + x_{13}^2}\\
0=& x_{8} + 16.44736842105263 (\cos(x_1)
x_{10} + \sin(x_1) x_{13} - x_3)\\
0=&  x_{9} + 10.319917440660475 (x_{4} - \sin(x_1)
x_{10} + \cos(x_1) x_{13})\\
0=&\sin(x_1) x_{8} + \cos(x_1) x_{9} - 0.030140727054618 (x_{10} -
 x_{11}) - 17.361008783459972 (x_{13} - x_{14})\\
0=&0.030140727054618 x_{10} - 1.395328440365198 x_{11} +
1.36518771331058 x_{12} + 17.361058783459974 x_{13} -\\ & \quad
28.877104346599904 x_{14} + 11.60409556313993 x_{15}\\
0=&1.36518771331058 (x_{11} -  x_{12}) +
11.60409556313993 x_{14} - 11.516095563139931 x_{15} - \frac{P
  x_{12}}{x_{12}^2 + x_{15}^2} - \frac{Q x_{15}}{x_{12}^2 +
  x_{15}^2}\\ 
0=&-(\cos(x_1) x_{8}) + \sin(x_1) x_{9} + 17.361008783459972 (x_{10} -
x_{11})-  0.030140727054618 (x_{13} -  x_{14})\\
0=&-17.361058783459974 x_{10} + 28.877104346599904 x_{11} -
11.60409556313993 x_{12} + 0.030140727054618 x_{13} \\ & \quad-
1.395328440365198 x_{14} + 1.36518771331058 x_{15}\\
0=&-11.60409556313993 x_{11} + 11.516095563139931 x_{12} +
1.36518771331058 (x_{14} - x_{15}) +
\frac{Q x_{12}}{x_{12}^2 + x_{15}^2} - \frac{P x_{15}}{x_{12}^2 +
  x_{15}^2} 
\end{align*}
\end{footnotesize}

The initial condition that gives a steady state is given by 
 \begin{footnotesize} $x(t_0)=[0.391057483977274, 376.9911184307751, \\ 1.022092319747551,
  0.308311065534821, 1.107019848098437, 0.199263572657719,
  1.12883036798339, 0.996801975949364,\\ 0.909203967958775,
  1.04, 1.006755413658047, 0.938198590465838,
  0, -0.070244002800643, -0.166824934470857]^T$\end{footnotesize}, here $P=1.25$ and
$Q=0.5$.

The stochastic noise is characterized by $\xi=[P,Q]^T$ and the
parameter sought is $\ipar$, while observing the voltage at the slack
bus, $x_{11}$ and $x_{14}$.

\end{document}